\documentclass{article}[12pt]
\pdfoutput=1
\usepackage[margin=1in]{geometry}
\usepackage[utf8]{inputenc}
\usepackage{authblk}
\usepackage{graphicx}
\usepackage{mathdots}
\usepackage{amsmath,amsthm,amssymb,mathrsfs}
\usepackage{centernot} %
\usepackage{stmaryrd}
\usepackage{graphicx}
\usepackage{tikz} %
\usetikzlibrary{shapes,decorations}
\usepackage{caption} %
\usepackage{subcaption} %
\usepackage{mdframed} %
\usepackage{upquote}
\usepackage{booktabs}
\usepackage{xspace}
\usepackage{amsbsy}
\usepackage{soul}
\usepackage[normalem]{ulem} %

\newcommand{\R}{\mathbb{R}}

\def\t{{^\text{T}}}

\begin{document}

\title{Learning partial differential equations for biological transport models from noisy spatiotemporal data}
\date{}
\author{%
John Lagergren$^{1,2,*}$, John T. Nardini$^{1,3,*}$, G. Michael Lavigne$^{1,2}$, Erica M. Rutter $^{1,2}$, Kevin B. Flores$^{1,2}$}

{\small{
\affil[1]{Department of Mathematics, North Carolina State University, Raleigh, North Carolina, United States of America}
\affil[2]{Center for Research in Scientific Computation, North Carolina State University, Raleigh, North Carolina, United States of America}
\affil[3]{The Statistical and Applied Mathematical Sciences Institute, Durham, North Carolina, United States of America}
\affil[*]{Contributed equally to this work}}}

\maketitle

\begin{abstract}
We investigate methods for learning partial differential equation (PDE) models from spatiotemporal data under biologically realistic levels and forms of noise. Recent progress in learning PDEs from data have used sparse regression to select candidate terms from a denoised set of data, including approximated partial derivatives. We analyze the performance in utilizing previous methods to denoise data for the task of discovering the governing system of partial differential equations (PDEs). We also develop a novel methodology that uses artificial neural networks (ANNs) to denoise data and approximate partial derivatives. We test the methodology on three PDE models for biological transport, i.e., the advection-diffusion,  classical Fisher-KPP, and nonlinear Fisher-KPP equations. We show that the ANN methodology outperforms previous denoising methods, including finite differences and polynomial regression splines, in the ability to accurately approximate partial derivatives and learn the correct PDE model. \\

\end{abstract}

\section{Introduction}

Recent research has investigated methods for discovering systems of differential equations that describe the underlying dynamics of spatiotemporal data. There are key advantages to learning and then using mathematical models for prediction instead of using a purely machine learning based method, e.g, neural networks. First, if the learned mathematical model is an accurate description of the processes governing the observed data, it has the ability to generalize from the set of training data to data outside of the training domain. Second, the learned mathematical model is interpretable, making it informative for scientists to hypothesize the underlying physical or biological laws governing the observed data. Examples of recent methods for inferring the underlying governing equations include the Sparse Identification of Nonlinear Dynamics (SINDy) algorithm \cite{brunton2016discovering} and the Equation Learner (EQL) neural network \cite{martius2016extrapolation,sahoo18learning}, both of which are used for discovering systems of ordinary differential equations (ODEs), and the PDE Functional Identification of Nonlinear Dynamics (PDE-FIND) algorithm \cite{rudy2017data}, which is used to identify PDE systems. Boninsegna et.~al \cite{boninsegna2018sparse} recently extended the SINDy algorithm to recover stochastic dynamical systems. Model selection criteria (such as Akaike Information Criteria and Bayesian Information Criteria) have been combined with the SINDy algorithm to increase robustness to errors, although incorrect models were still selected at noise levels we consider here \cite{mangan2017model}. 
The discovery methods mentioned above assume that the measured data arise from a parameterized $n$-dimensional dynamical system of the form
\begin{subequations}
    \begin{align}
        u_t(x,t) &= F\left(x,t,u,u_x,u_{xx},\dots,\theta\right), & x\in[x_0,x_f], \quad t\in[t_0,t_f] \label{math_model}\\ 
        u(x,t_0) &= u_0(x), & x\in[x_0,x_f] \label{init_cond}
    \end{align}
\end{subequations}
with parameter vector $\theta\in\R^{k}$ and appropriate boundary conditions. It is then assumed that the true dynamical system, $F(\cdot)$, is only comprised of a few terms (e.g., for the diffusion-advection equation, $F=Du_{xx}-cu_{x},D,c\in\R$). The goal of these methods is to correctly specify the small number of correct terms from a large library with many potential candidate terms (such as $u$, $u_x$, $u_{xx}$, $u^2$, $uu_x$, etc.) with the aid of sparse regression \cite{Rish_sparse_2015}. A practical challenge arises because one typically does not have access to the noiseless values of $u(x,t)$ or its partial derivatives. Instead, one must approximate these values from noisy experimental data. Here, our goal is to investigate the performance of existing denoising methods that are used in conjunction with PDE-FIND and to present a novel denoising methodology relying on artificial neural networks (ANNs).

Several methods have been used for denoising data to approximate $u(x,t)$ and it partial derivatives ($u_t$, $u_x$, $u_{xx}$, etc.). The most prevalent methods that have been proposed are finite difference approximations or the use of cubic splines for interpolation, followed by partial differentiation of the fitted splines. However, both of these methods have been found to be prone to inaccuracies in the presence of noise \cite{rudy2017data}. Recent work has considered recovery of dynamical systems with high amounts of noise added to the time derivative measurement ($u_t$) by transforming the data into a spectral domain \cite{schaeffer2017learning}. Zhang et al.~\cite{zhang2018robust} recently proposed using sparse Bayesian regression, which allows for error bars for each candidate term in the discovered equation. However, although their method was robust, the noise in this study was also added only to the time derivative term ($u_t$) instead of the observed data ($u(x,t)$). %
Importantly, it has been noted that introducing noise to the observed data itself ($u(x,t)$) hinders the recovery of the correct PDE, thus, developing a method of denoising data for $u(x,t)$ has been identified as a current challenge for learning PDEs \cite{schaeffer2017learning}. To the best of our knowledge, the use of finite differences or utilizing splines are the two methods that, in practice, yield the most accurate approximations for the library terms involved in PDE learning. The primary challenge involved with using these methods for numerical differentiation, which we further test in this work, is that they are sensitive to noise levels and can amplify noise as the order of the derivative increases. This challenge inhibits learning PDEs for practical biological applications where data may have large noise levels due to many sources of error, including the data collection process, imprecise measurement tools, and the inherent stochastic nature of biological processes \cite{banks2014modeling,codling_edward_a_random_2008}. For example, for ecological measurements of population abundance, typical data sets can have noise levels on the order of a coefficient of variation equal to 0.2 \cite{francis_quantifying_2003,perretti_model-free_2013}.
Notably, adding this biologically relevant level of noise to the observation $u(x,t)$ has not been considered in previous PDE learning work.

An additional, yet realistic, complication that has not been considered is the presence of non-constant error noise in the spatiotemporal data used for PDE learning. For example, proportional error noise can occur when the variance of the data is proportional to the size of the measurement, e.g., population size or density \cite{banks_estimation_2011}. Non-constant error noise may also occur when the observed processes occur on different time scales \cite{johnston_how_2014}. To account for non-constant error noise in the scenario that one has a mathematical model for the biological process generating the data, e.g., $u(x,t)$, the non-constant error noise can be accounted for with a statistical model used in conjunction with the mathematical model \cite{anderssen_numerical_1974}. For example, for a set of observed data at space points $x_i$, $i=1,\dots,M$ and time points $t_j$, $j=1,\dots,N$, a general statistical model is given by 
\begin{align}
    U_{i,j} &= u(x_i,t_j) + w_{i,j}\odot\mathcal{E}_{i,j}, \label{stat_model}
\end{align}
where the noiseless observations are corrupted by noise modeled by the random variable $w_{i,j}\odot\mathcal{E}_{i,j}$ in which $\odot$ represents element-wise multiplication. Finite difference methods assume $w_{i,j}\odot\mathcal{E}_{i,j}=0$ while regression methods using splines assume the variance of $w_{i,j}\odot\mathcal{E}_{i,j}$ is constant. More generally, the error term $\mathcal{E}_{i,j}$ may instead be generated by a probability distribution that is weighted by 
\begin{equation}
    w_{i,j} = \left(\beta_1u^\gamma_1(x_i,t_j), \dots, \beta_nu^\gamma_n(x_i,t_j) \right)\t \label{stat_model2}
\end{equation}
for $\gamma\geq0$ and $\beta_1,\dots,\beta_n\in\R$. Constant error noise is modeled by assuming $\gamma=0$ and $\beta_1,\dots,\beta_n=1$. Proportional error noise is modeled by assuming $\gamma>0$, $\beta_1,\dots,\beta_n\neq0$ \cite{banks_mathematical_2009}. We note that finite difference and spline approximations of noisy data neglect whether the noise process has non-constant variance. We hypothesize that this leads finite difference and spline approximations to yield poor estimates of the noiseless data $u(x,t)$ and its partial derivatives when the data contain proportional error noise. In this work, we investigate this hypothesis and develop a methodology using ANNs as a model for $u(x,t)$ in conjunction with an appropriate statistical error model that accounts for the presence of proportional error when denoising spatiotemporal data.

The denoising methods present in this work focus on spatiotemporal data for learning PDEs, however the methods we describe can be readily applied to learning ODEs. We choose to focus our study on a specific set of diffusive PDE models, which have provided a wealth of insight into many biological transport phenomena, including ecological migration and invasion \cite{hastings_spatial_2005}, neuronal transport \cite{jones_stability_1984}, cancer progression \cite{baldock_patient-specific_2014,kuang_data-motivated_2015,rockne_patient-specific_2015,rutter_mathematical_2017}, and wound healing \cite{maini_travelling_2004,nardini_modeling_2016}. 
Here, we demonstrate how an ANN can be used with a non-constant error statistical model to accurately approximate $u(x,t)$ from noisy proportional error data.
It has long been known that ANNs are universal function approximators \cite{hornik_1991}, meaning ANNs have the capacity to approximate continuous functions arbitrarily well. Unlike local approximations such as finite differences and splines, an ANN can be fit to an entire set of spatiotemporal data, affording a global context that may help to decrease overfitting. We compare the accuracy of the ANN-based method to finite differences and splines in computing estimates for $u(x,t)$ and its partial derivatives. We then investigate whether the PDE-FIND algorithm can more accurately infer the governing PDE equations from data when its library of terms is constructed using an ANN-based method.

\begin{figure}
    \centering
    \includegraphics[width=0.89\textwidth]{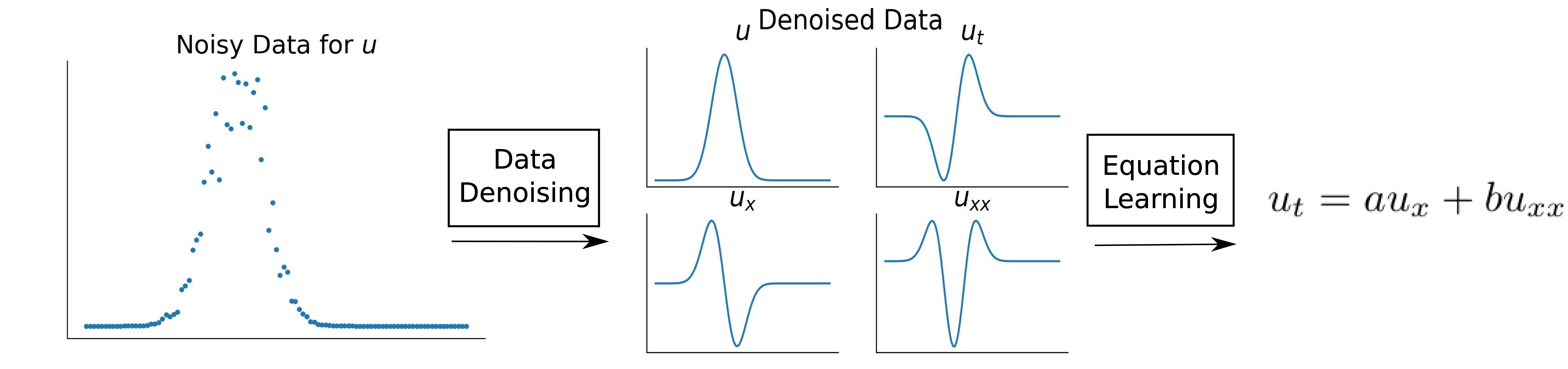}
    \caption{The two components to learning PDEs from data.  The first component is to approximate $u$, $u_t$, $u_x$, $u_{xx}$, etc., from noisy data. The second component uses the output from the first component as an input for the PDE-FIND algorithm to learn a PDE. We also employ a pruning algorithm after the PDE-FIND step, not depicted here.}
    \label{pipeline}
\end{figure}

\section{Methods} \label{Methods}

The process of learning a system of equations from noisy data can be divided into two main components: (1) the data denoising and library construction component, in which the underlying dynamical system $u(x,t)$ and its partial derivatives are approximated from the noisy realizations of $\big\{U_{i,j}\big\}_{i=1,j=1}^{M,N}$ in \eqref{stat_model}, and (2) the equation learning component, in which, given approximations for $u$, $u_t$, $u_x$, $u_{xx}$, etc., one employs an algorithm that can effectively uncover the mechanistic form of $F$ in \eqref{math_model} (Figure \ref{pipeline}). Below, we describe the mathematical models used for data generation, the method of constructing a library from noisy data, and equation learning. All of the denoising methods were implemented in Python 2.7 using the Scipy package for polynomial splines and the Keras machine learning library for ANNs. All code and videos are available at \verb+https://github.com/biomathlab/PDElearning/+.

\subsection{Data Generation}
We consider three diffusive PDE models for biological transport in this work, each of which has been used previously to interpret biological data \cite{jin_reproducibility_2016,maini_travelling_2004,sherratt_models_1990,sibert_advectiondiffusionreaction_1999}. These models include the diffusion-advection equation
\begin{equation}
u_t = - cu_x + Du_{xx}, \hspace{.5cm}D,c\in\R \label{diffusion_advection}
\end{equation}
the classical Fisher-Kolmogorov-Petrovsky-Piskunov (Fisher-KPP) Equation
\begin{equation}
    u_t=Du_{xx}+ru - ru^2,\hspace{.5cm}D,r\in\R \label{fisher_kpp}
\end{equation}
and the nonlinear Fisher-KPP Equation
\begin{equation}
    u_t= Duu_{xx} + Du_x^2 + ru -ru^2,\hspace{.5cm}D,r\in\R\label{nonlinear_fisher}
\end{equation}
where $D$ is the diffusion coefficient, $r$ is the intrinsic population growth rate, and $c$ is the advection rate.

Assume $u(x,t)$ denotes the solution to one of the above mathematical models. We generate noisy data by using Equation \eqref{stat_model} with $w_{i,j}=\sigma u(x_i,t_j)$ (i.e., $\beta = \sigma \text{ and } \gamma=1$) in which all $\mathcal{E}_{i,j}$ terms are simulated as i.i.d. normal random variables with mean zero and variance one. We generate six data sets for each mathematical model, setting $\sigma=$ 0, 0.01, 0.05, 0.10, 0.25, and 0.50. For Equation \eqref{diffusion_advection}, we use its analytical solution to compute $u(x,t)$. For Equations \eqref{fisher_kpp} and \eqref{nonlinear_fisher} we use finite difference computations to numerically approximate $u(x,t)$. The numerical step sizes in these computations were chosen small enough to not introduce significant noise into the solution. We use $M=101$ spatial locations and $N = 300$ time points to generate data for the diffusion-advection equation and $M=199$ spatial locations and $N=99$ time points for the Fisher-KPP and nonlinear Fisher-KPP Equations. %
As a preprocessing step, each data set is scaled to $[0,1]$ in order to consistently measure errors across the various data sets and noise levels. 

\subsection{Data denoising and library construction}

\subsubsection{Finite differences and spline approximations}

For finite differences, we use central difference formulas on interior points and forward differences at the boundaries to obtain first order derivative approximations. For higher order derivatives (e.g. $u_{xx}$), first order finite difference rules are repeated on the corresponding previous order derivative approximations. 

Finite difference approximations can be obtained accurately, efficiently, and directly from noiseless data, however, their accuracy quickly deteriorates in the presence of observation error. Following \cite{rudy2017data}, we also employ polynomial spline regression. For a given data point $u_{i,j}$ in the set of observations, we fit a cubic bi-spline on a small two-dimensional neighborhood of size $11\times11$ centered at $u_{i,j}$. Denoised function and derivative approximations are then obtained using evaluations and the analytic derivatives of the fitted polynomial at the center point $u_{i,j}$. Note that since we only approximate derivatives up to second order, higher order splines are not considered. For points close to the boundaries, we evaluate along the spline approximation that is nearest to the boundary. Each polynomial spline is optimized using ordinary least squares. We found that cubic bi-spline approximations were generally more robust than the one-dimensional cubic splines used in \cite{rudy2017data} (Supplementary Figures S1-S3). 

\subsubsection{Artificial neural network approximations}

An artificial neural network (ANN), denoted $h_\theta(\vec{x})$, was used to approximate $u(x,t)$, with one hidden layer of the form
\begin{equation}
    h_\theta(\vec{x}) = a_2\bigg(W_2\Big(a_1(W_1\vec{x} + b_1)\Big) + b_2\bigg) \label{ann}
\end{equation}
where $\vec{x}=[x\,\,\,t]\t$ and $a_i(\cdot)$ represents the continuous nonlinear activation function for the i$^\text{th}$ layer. The matrices $W_i$ and vectors $b_i$ (typically called weights and biases) comprise the total set of trainable parameters $\theta=\{W_1,b_1,W_2,b_2\}$ of the network. We note that the use of one hidden layer in a neural network is a sufficient condition to make it a universal function approximator under the assumption that the activation function is bounded and non-constant \cite{hornik_1991}. This result extends to ANNs with multiple hidden layers, however, we found that while training multi-layer ANNs resulted in faster convergence, the derivative approximations were worse. The task is therefore to find the optimal parameters $\theta^*$ such that $h_{\theta^*}(x,t) \approx u(x,t)$. The fitted surface function $h_{\theta^*}(x,t)$ and the computation of analytic derivatives of this function are used to approximate $u(x,t)$ in \eqref{stat_model} and its partial derivatives for library construction in the PDE learning task.

We formulate a regression problem using the ANN defined in \eqref{ann} in order to find $\theta^*$ that minimizes the generalized least squares cost function
\begin{equation}
    \mathcal{J}(\theta) = \frac{1}{MN}\sum_{i=1,j=1}^{M,N} \left( \frac{h_\theta(x_i,t_j) - u_{i,j}}{|h_\theta(x_i,t_j)|^\gamma} \right)^2. \label{cost}
\end{equation}
Note that this error formulation accounts for the statistical error model in \eqref{stat_model} and also reduces to ordinary least squares when $\gamma=0$. Also note that function values $h_\theta(x_i,t_j)$ less than 1e-4 in absolute value are set equal to one in the denominator $|h_\theta(x_i,t_j)|^\gamma$ during network training and evaluation. The network parameters $\theta$ are optimized using the first-order gradient-based ``Adam'' optimizer \cite{kingma2014adam} with default parameters and a batch size of 10.

We found that the choice of activation function in the ANN, $a_i(\cdot)$, plays an important role in the accuracy of the partial derivative approximations. Typical activations like sigmoid and hyperbolic tangent yield oscillations in higher order derivative terms (\emph{e.g}., see Supplementary movie S1.). To mitigate this, we chose to use the ``softplus'' activation function which takes the form $\log(1+e^z)$. This function has many desirable properties for approximating $u(x,t)$ (e.g. smoothness and infinitely many derivatives), however, it is unbounded, which violates an assumption of ANNs as universal approximators. We address the unboundedness of the softplus function by including an $\ell_2$-regularization penalty on the inputs $W_1\vec{x} + b_1$ in \eqref{ann}. Further, an additional squared error term is included in the loss function to penalize function values outside $[0,1]$. Without this term, the function values can blow up during training since $\mathcal{J}(\theta)\rightarrow1$ as $h_\theta(\vec{x})\rightarrow+\infty$ when $\gamma=1$. Thus, the complete loss function used for training the ANN is
\begin{equation}
\mathcal{L}(\theta) = \frac{1}{MN} \sum_{i=1,j=1}^{M,N} \left[\left( \frac{h_\theta(x_i,t_j) - u_{i,j}}{|h_\theta(x_i,t_j)|^\gamma} \right)^2 + \lambda \left(W_1\vec{x}_{i,j} + b_1\right)^2 \right] + \frac{1}{MN}\sum_{h_\theta\not\in[0,1]} h_\theta^2  \label{loss}
\end{equation}
where the first term corresponds to the generalized least squares cost function \eqref{cost}, the second term corresponds to the $\ell_2$-regularization penalty with regularization strength $\lambda$, and the third term corresponds to the additional squared error term to penalize function values outside $[0,1]$. Generally, the regularization strength should be chosen separately for each data set using cross-validation, but we found that no regularization (i.e. $\lambda=0$) was needed for the 18 data sets considered in this paper. 

We used 1,000 neurons in the hidden layer of the ANN. This choice was large enough to have maximal capacity to fit the data, while still allowing the optimization of $\theta$ to be computationally feasible on a desktop computer (3.4 GHz Intel Core i5 processor, 8gb RAM). In order to prevent overfitting, the data were randomly split into training and validation sets (90\%/10\%) when training the ANN. The optimal network parameters were chosen to be those that minimized the error \eqref{cost} on the validation set. Training ended if the validation error did not decrease for 50 consecutive epochs. 

Representative examples of results from the bi-spline and ANN methods are shown in Supplementary movies S2 and S3, respectively. All movies for all methods and noise levels considered in this work can be found at \verb+https://github.com/biomathlab/PDElearning/animations/+.

\subsection{Equation Learning}\label{eql_methods}

We use the PDE-FIND algorithm \cite{rudy2017data} to discover the form of $F(\cdot)$ in Equation \eqref{math_model} using computations of $u,u_x,u_{xx},$ and $u_t$ from the ANN, spline, and finite difference methods. Prior to implementing PDE-FIND, the numerical approximations are scaled from $[0,1]$ back into their original scales. We discuss the PDE-FIND implementation in Section \ref{PDEFIND implementation} and an additional pruning method in Section \ref{pruning} that is used to remove extra terms from the final learned equation. We further discuss how we analyze our results in Section \ref{implementation}.

\subsubsection{PDE-FIND Implementation} \label{PDEFIND implementation}

Once $u(x,t)$ and its partial derivatives have been computed,  a large library of potential PDE terms is formed column-wise in the matrix, $\Theta$, given by

\begin{equation}
    \Theta = 
    \begin{bmatrix}
        1 & u & \cdots & u^p & u_x &  \cdots & u^p\odot u_x & u_{xx} & \cdots & u^p\odot u_{xx} & u_x^2 & u_x\odot u_{xx} & u_{xx}^2
    \end{bmatrix}
\end{equation}
where each column of $\Theta$ is some vectorization of the written term. The ANN has difficulty capturing the early dynamics of the diffusion-advection equation, so we skip the first 20 timepoints from the denoised data when building $\Theta$ for all data sets and denoising strategies. To reduce the computational time, only every fifth remaining timepoint is included in $\Theta$. Hence, while the data sets for the diffusion-advection equation begin with $N=300$ timepoints, only (300-20)/5=56 timepoints are used in constructing $\Theta$. We set $p=2$ resulting in $d=12$ columns in $\Theta$. %
Each column of $\Theta$ thus represents a candidate term comprising $F$, so we assume 

\begin{equation}
    u_t \approx \Theta \xi, \label{linear_regression}
\end{equation}
where $\xi$ is a vector whose nonzero entries correspond to the true terms of $F$. The vector $\xi$ is estimated using methods from sparse regression \cite{Rish_sparse_2015}. Sequential threshold ridge regression was found to be a suitable method for estimating $\xi$ for PDE-FIND in a previous study \cite{rudy2017data}. However, we found that the Greedy algorithm to performed well for the data and models we considered in this work. The Greedy algorithm computes

\begin{equation}
    \hat{\xi} = \arg\min_{\xi\in\mathbb{R}^d} \dfrac{1}{MN} \|u_t - \Theta\xi\|_2^2, \quad \text{subject to } \|\xi\|_0\leq k
\end{equation}
for some sparsity parameter, $k$ \cite{zhang2009adaptive}.

The value of $k$ for a given data set is treated as a hyperparameter that is found by splitting the library data into separate training and validation sets, and then optimizing over the validation set. In this training-validation procedure, we randomly divide our data points for $u_t$ into 5-by-5 tiles of adjacent spatiotemporal points and then randomly assign 50\% of these tiles to a training data set, $u_t^{\text{train}}$, and the remaining 50\% to a validation set, $u_t^{\text{validate}}$. We split the corresponding rows of $\Theta$ into $\Theta^{\text{train}}$ and $\Theta^{\text{validate}}$. We perform our hyperparameter search over 51 values for $k$ between 0 and $10^3$. For each value of $k$, we estimate $\hat{\xi}$ from the training set. For each estimate, we then compute its mean-squared error (MSE) over the validation set. We choose the hyperparameter corresponding to the $\hat{\xi}$ estimate with the smallest MSE on the validation data. The equation that results from sparse regression with this hyperparameter is our final equation from the PDE-FIND algorithm. We refer to the validation MSE from the final equation ``$\text{val}_0$'' in the remaining text.

\subsubsection{Pruning Method} \label{pruning}

We chose a 50-50 training and validation split for the data to avoid overfitting to the training data with a large validation set. Even so, we will demonstrate in Section \ref{PDE-FIND_pruning} below that PDE-FIND is able to learn small but systematic biases from the ANN's fit to $u$ and its derivatives by incorporating extra terms into the final equations. Pruning methods have previously been developed that remove extra terms that do not significantly increase an algorithm's performance, see for example \cite{anders_model_1999,mackay_probable_1995}.  Accordingly, we implement the following pruning method after the PDE-FIND implementation described in Section \ref{PDEFIND implementation} for all methods in order to delete the extra terms from the final equation.  

The pruning procedure starts with the reduced library of candidate terms (i.e., the columns of $\Theta$)  for the right hand side of Equation (\ref{math_model}) with nonzero entries of $\hat{\xi}$ that resulted from our training-validation procedure. We then perform a sensitivity test for the remaining terms as follows. Suppose $\tilde{d}$ terms remain in $\tilde{\Theta}$, and let $\tilde{\Theta}_i, i=1,\hdots,\tilde{d}$ denote the further-reduced library where the $i$th column of $\tilde{\Theta}$ has been removed. For each value of $i$, we find the least squares solution (without regularization) on the training data to the equation

\begin{equation}
    u_t^{\text{train}} = \tilde{\Theta}_i^{\text{train}}\xi_i.
\end{equation}
We then use our $\hat{\xi}_i$ estimate and  compute the MSE over the validation data when the $i^\text{th}$ term has been removed and call this computation $\text{val}_i$. We then remove any candidate terms for our library that result in $\text{val}_i/\text{val}_0<1+\alpha$ for some $\alpha>0$. After this pruning step, we perform one final round of training without regularization over the fully reduced library to find the final form of our underlying equation.

It is important to note that choice of the $\alpha$ pruning threshold value warrants careful decision. If this value is chosen too high, then too few terms will be selected and the learned equation will be incomplete. If the chosen value is too small, then the final equation will admit extra terms arising from the systematic errors in derivative estimation. We will demonstrate below that the choice of arbitrary choice of $\alpha = 0.25$ below provides promising results for the diffusion-advection and Fisher-KPP Equations, while $\alpha=0.05$ is suitable for the nonlinear Fisher-KPP Equation.

\subsubsection{Accuracy Metrics} \label{implementation}

To quantitatively assess the accuracy in recovering the correct PDE that generated the data, i.e., using the combined PDE-FIND with pruning methodology described above, we introduce the the \emph{true positive ratio} (TPR) for a given vector $\xi$ as:

\begin{equation}
    \text{TPR}(\xi) = \dfrac{\text{TP}}{\text{TP + FN + FP}}, \label{TPR}
\end{equation}
where TP (``True Positive'') denotes the number of correctly-specified nonzero coefficients in $\xi$, FN (``False Negative'') denotes the number of coefficients in $\xi$ that are incorrectly specified as zero, and FP (``False Positive'') denotes the number of coefficients in $\xi$ that are incorrectly specified as nonzero. Recall that the nonzero entries of $\xi$ correspond to the relevant terms in an equation (i.e., for a library of $\Theta=[1 \  u \  u_x \  u_{xx}], \xi = [0 \  1 \  2 \  0]^T$ corresponds to the equation $u_t = u + 2u_x$). For example, when trying to learn Equation \eqref{diffusion_advection}, an equation of the form $u_t = u_{xx} + uu_x$ would have TP = 1 (the nonzero coefficient for $u_{xx}$ is correct), FN = 1 (the missing $u_x$ term is incorrect), FP = 1 (the nonzero $uu_x$ term is incorrect), resulting in a final score of TPR = 1/3. Note that the TPR value is similar in nature to the Jaccard index: larger TPR values suggest that the true equation form has been better approximated, and TPR = 1 signifies that the correct equation form has been recovered. 

We note that the learned equation from the PDE-FIND with pruning method was often found to be sensitive to the random split of $u_t$ and $\Theta$ into training and validation data. Therefore, we performed PDE learning for 1,000 different random training-validation data splits of $u_t$ and $\Theta$ for each data set and for each computational method (finite differences, splines, and ANN). We then consider the distribution of TPR($\hat{\xi}$) scores to assess the overall performance of the methodology. We declare the most commonly-learned equation among the 1,000 data splits as the final learned equation for each data set and computational method.
\begin{table}[ht]
    \centering
    \caption{The relative mean-squared error (RMSE) between the noiseless data or true derivative values and our desnoised data or derivative computations using finite differences, splines, and the ANN for the diffusion-advection equation. ``FD'' denotes finite differences, ``SP'' denotes splines, and ``ANN'' denotes the ANN. Bold denotes the lowest errors of the three methods.}
    \label{advection_diffusion_rmse}
    \begin{tabular}{cccccc}
        Noise & Method & $u$ RMSE & $u_{t}$ RMSE & $u_{x}$ RMSE & $u_{xx}$ RMSE \\ 
        \hline
         & FD & \textbf{0.00e+00} & \textbf{5.39e-05} & \textbf{5.77e-04} & \textbf{3.69e-02} \\
        $\sigma$=0 & SP & 1.22e-02 & 1.12e+00 & 3.33e-02 & 4.86e+01 \\
         & ANN & 2.86e-04 & 3.96e-01 & 8.47e-03 & 3.75e-01 \\
        \hline
         & FD & \textbf{1.02e-04} & 2.08e+02 & 4.34e-01 & 3.52e+01 \\
        $\sigma$=0.01 & SP & 1.11e-02 & 7.34e+00 & 4.93e-02 & 4.87e+01 \\
         & ANN & 8.40e-04 & \textbf{7.71e-02} & \textbf{1.15e-02} & \textbf{6.93e-01} \\
        \hline
         & FD & 2.51e-03 & 2.42e+03 & 9.97e+00 & 1.01e+03 \\
        $\sigma$=0.05 & SP & 1.19e-02 & 7.47e+01 & 2.62e-01 & 5.01e+01 \\
         & ANN & \textbf{5.61e-04} & \textbf{1.95e-01} & \textbf{7.71e-03} & \textbf{7.90e-01} \\
        \hline
         & FD & 1.00e-02 & 4.04e+03 & 7.59e+01 & 3.78e+03 \\
         $\sigma$=0.10 & SP & 1.04e-02 & 2.32e+02 & 1.38e+00 & 5.45e+01 \\
         & ANN & \textbf{9.51e-04} & \textbf{1.23e-01} & \textbf{1.44e-02} & \textbf{7.69e-01} \\
        \hline
         & FD & 6.28e-02 & 9.04e+04 & 2.20e+02 & 3.65e+04 \\
         $\sigma$=0.25 & SP & 2.51e-02 & 6.21e+03 & 5.83e+00 & 2.42e+02 \\
         & ANN & \textbf{7.29e-03} & \textbf{1.53e+00} & \textbf{4.49e-02} & \textbf{7.21e-01} \\
        \hline
         & FD & 2.41e-01 & 2.58e+06 & 1.39e+03 & 1.05e+05 \\
         $\sigma$=0.50 & SP & \textbf{3.71e-02} & 3.78e+04 & 6.68e+01 & 7.95e+02 \\
         & ANN & 6.34e-02 & \textbf{3.43e+00} & \textbf{1.05e-01} & \textbf{1.44e+00} \\
        \hline
    \end{tabular}
\end{table}

\section{Results} \label{Results}

In this section, we detail our results using the ANN to denoise data for $u(x,t)$ and compute partial derivatives. In addition, we test the accuracy of using the ANN method in conjunction with PDE-FIND to learn PDEs. Analogous results are presented for finite differences and splines. We begin by demonstrating the accuracy of the partial derivative calculations in Section \ref{Derivative_calc}, we explain why PDE-FIND finds small systematic bias terms in Section \ref{PDE-FIND_pruning}, and then we detail the accuracy in learning of the diffusion-advection, Fisher-KPP, and nonlinear Fisher-KPP equations in Sections \ref{PDE_find_DA}-\ref{PDE_find_nonlinear_FKPP}.

\subsection{Derivative Calculations}\label{Derivative_calc}

\begin{table}[ht]
    \centering
    \caption{Relative mean square error for Fisher-KPP equation.   ``FD'' denotes finite differences, ``SP'' denotes splines, and ``NN'' denotes the ANN. Bold denotes the lowest errors of the three methods.}
    \label{fisher_rmse}
    \begin{tabular}{cccccc}
        Noise & Method & $u$ RMSE & $u_{t}$ RMSE & $u_{x}$ RMSE & $u_{xx}$ RMSE \\ 
        \hline
         & FD & \textbf{0.00e+00} & \textbf{3.16e-05} & \textbf{4.46e-04} & \textbf{4.67e-03} \\
         $\sigma=0$ & SP & 6.28e-05 & 2.83e-04 & 2.83e-03 & 2.01e-01 \\
         & ANN & 4.86e-04 & 6.98e-02 & 1.18e-01 & 2.66e+00 \\
        \hline
         & FD & 9.82e-05 & 4.84e+00 & 6.47e+01 & 1.19e+03 \\
        $\sigma=0.01$ & SP & \textbf{6.91e-05} & 1.89e-01 & 1.56e+00 & \textbf{2.19e+00} \\
         & ANN & 3.90e-04 & \textbf{3.85e-02} & \textbf{1.19e-02} & 2.55e+00 \\
        \hline
         & FD & 2.52e-03 & 9.49e+01 & 9.43e+02 & 6.78e+04 \\
         $\sigma=0.05$ & SP & \textbf{2.20e-04} & 8.85e+00 & 5.83e+01 & 2.39e+02 \\
         & ANN & 4.67e-04 & \textbf{1.03e-02} & \textbf{1.63e-02} & \textbf{1.57e+00} \\
        \hline
         & FD & 9.95e-03 & 4.05e+02 & 2.81e+03 & 2.02e+05 \\
         $\sigma=0.10$ & SP & \textbf{6.80e-04} & 1.64e+01 & 6.66e+01 & 9.35e+02 \\
         & ANN & 8.46e-04 & \textbf{4.26e-02} & \textbf{5.29e-02} & \textbf{1.23e+00} \\
        \hline
         & FD & 6.30e-02 & 2.92e+03 & 3.29e+04 & 3.92e+05 \\
         $\sigma$=0.25 & SP & \textbf{4.14e-03} & 2.26e+02 & 5.65e+02 & 1.44e+04 \\
         & ANN & 6.41e-03 & \textbf{6.94e-02} & \textbf{1.04e-01} & \textbf{5.90e+00} \\
        \hline
         & FD & 2.38e-01 & 9.22e+03 & 1.01e+05 & 1.04e+06 \\
         $\sigma=0.50$ & SP & \textbf{1.52e-02} & 5.44e+02 & 1.32e+03 & 3.22e+04 \\
         & ANN & 6.60e-02 & \textbf{5.48e-01} & \textbf{8.98e-01} & \textbf{3.52e+01} \\
        \hline
    \end{tabular}
\end{table}

\begin{table}[ht]
    \centering
    \caption{Relative mean square error for nonlinear Fisher-KPP equation.   ``FD'' denotes finite differences, ``SP'' denotes splines, and ``ANN'' denotes the ANN. Bold denotes the lowest errors of the three methods.}
    \label{fisher_nonlin_rmse}
    \begin{tabular}{cccccc}
        Noise & Method & $u$ RMSE & $u_{t}$ RMSE & $u_{x}$ RMSE & $u_{xx}$ RMSE \\ 
        \hline
         & FD & \textbf{9.71e-36} & \textbf{2.20e-05} & \textbf{1.53e-04} & \textbf{1.34e-01} \\
         $\sigma=0$ & SP & 1.53e-05 & 4.19e-05 & 1.09e-03 & 5.71e+00 \\
         & ANN & 8.73e-04 & 1.41e+01 & 6.16e+00 & 1.40e+02 \\
        \hline
         & FD & 1.02e-04 & 2.21e+02 & 9.16e+03 & 8.26e+02 \\
        $\sigma=0.01$ & SP & \textbf{2.01e-05} & 1.02e+01 & 6.20e+02 & \textbf{3.32e+00} \\
         & ANN & 6.87e-04 & \textbf{5.55e+00} & \textbf{6.58e+01} & 1.90e+02 \\
        \hline
         & FD & 2.43e-03 & 5.45e+03 & 2.65e+05 & 3.35e+04 \\
         $\sigma=0.05$ & SP & \textbf{1.57e-04} & 2.74e+02 & 6.84e+03 & \textbf{6.69e+01} \\
         & ANN & 1.08e-03 & \textbf{8.26e+00} & \textbf{1.68e+00} & 1.70e+02 \\
        \hline
         & FD & 1.01e-02 & 2.30e+04 & 1.16e+06 & 5.69e+04 \\
         $\sigma=0.10$ & SP & \textbf{6.20e-04} & 1.32e+03 & 1.73e+04 & \textbf{7.23e+01} \\
         & ANN & 1.84e-03 & \textbf{1.93e+01} & \textbf{1.10e+01} & 2.04e+02 \\
        \hline
         & FD & 6.25e-02 & 1.47e+05 & 4.95e+06 & 9.86e+05 \\
         $\sigma=0.25$ & SP & \textbf{4.00e-03} & 5.69e+03 & 1.24e+05 & 1.54e+03 \\
         & ANN & 6.34e-03 & \textbf{2.58e+01} & \textbf{2.34e+01} & \textbf{2.59e+02} \\
        \hline
         & FD & 2.43e-01 & 5.73e+05 & 1.92e+07 & 2.46e+06 \\
        $\sigma=0.50$ & SP & \textbf{1.38e-02} & 2.87e+04 & 4.16e+05 & 1.13e+04 \\
         & ANN & 6.89e-02 & \textbf{6.49e+01} & \textbf{1.97e+02} & \textbf{4.88e+02} \\
        \hline
    \end{tabular}
\end{table}

We found that the finite difference method most accurately approximates $u$ and its derivatives for $\sigma=0$ (Table \ref{advection_diffusion_rmse}). This result is not surprising, as finite difference computations assume that there is no error in the data. For all other values of $\sigma$, we observe that the ANN produces the most accurate derivative calculations, though the splines appear best at inferring $u(x,t)$ when the data are noisy. It is important to note that the ANN's calculations are often several magnitudes of order more accurate than the spline and finite difference approximations, and this disparity between the computations appears to increase with $\sigma$ (Table \ref{advection_diffusion_rmse}). For example, at $\sigma=0.01$, the ANN's relative mean squared error (RMSE) for $u_t$ is four orders of magnitude smaller than the RMSE for finite differences and two orders of magnitude smaller than the RMSE for splines. At $\sigma=0.50$, the ANN's RMSE for $u_t$ has become six orders of magnitude smaller than the RMSE for finite differences and four orders of magnitude smaller than the RMSE for splines. The other derivative computations show similar results.

Similarly, we found that the ANN was most accurate for noisy data from the Fisher-KPP Equation (Table \ref{fisher_rmse}) and the nonlinear Fisher-KPP Equation (Table \ref{fisher_nonlin_rmse}). Recall that we do not have analytical solutions to these equations, so we used finite difference computations on the noiseless data $(\sigma=0)$ as an estimate for the analytical derivative values for the Fisher-KPP and nonlinear Fisher-KPP Equations. Again, in both cases we observe that the finite difference calculations perform best in computing the RMSE for $\sigma=0$, but on average, the ANN provides the best calculations for the derivatives for larger values of $\sigma$. The spline method is consistently the most accurate at inferring $u$ from the data. The disparity between the RMSE calculations for the ANN as compared to the splines or finite differences again appears to increase with $\sigma$ for these two equations.

\subsection{PDE-FIND without pruning learns the wrong equation} \label{PDE-FIND_pruning}
\begin{figure}[h!]
    \centering
    \includegraphics[width=0.4\textwidth]{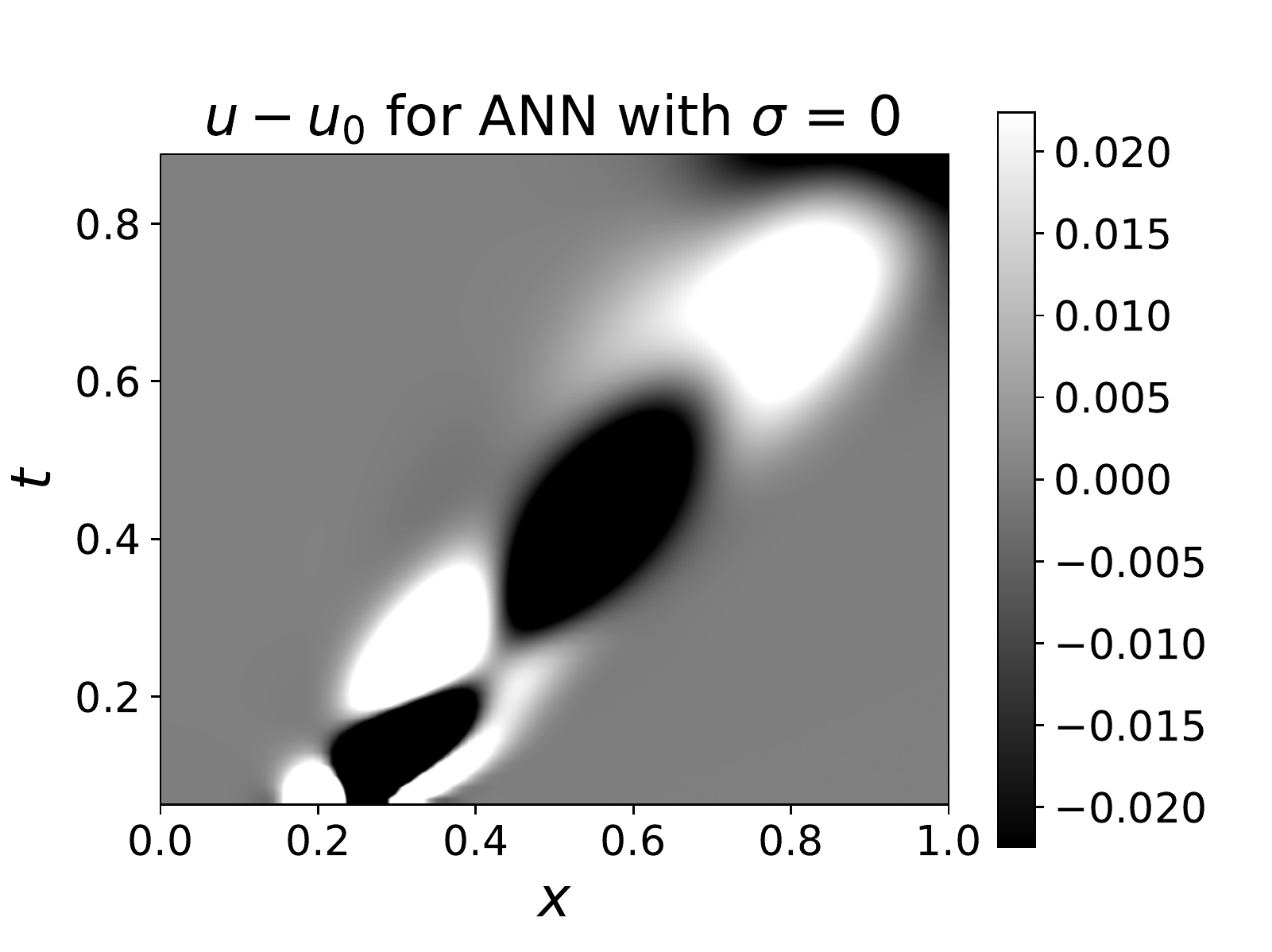}\includegraphics[width=0.4\textwidth]{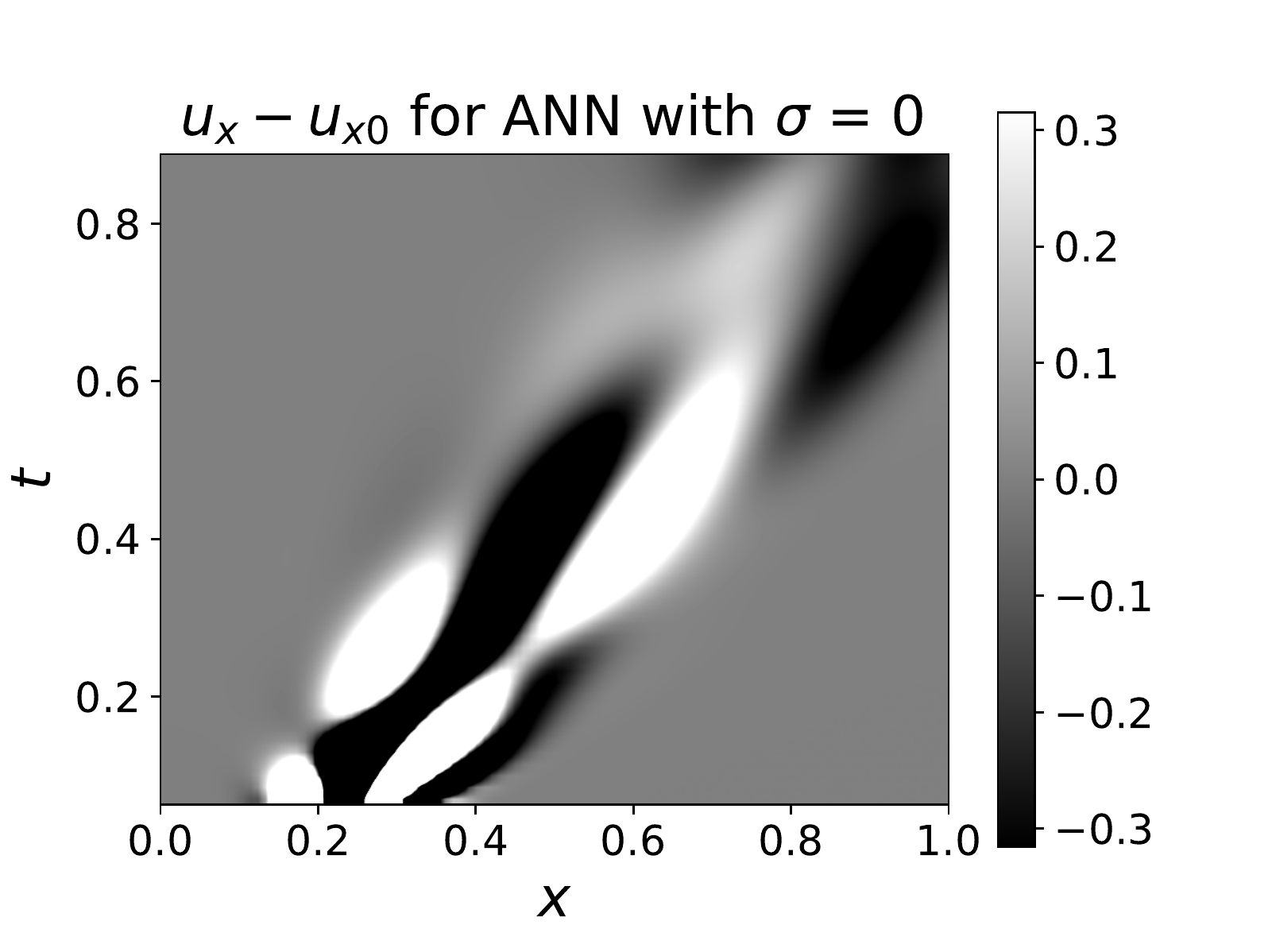}
    \includegraphics[width=0.4\textwidth]{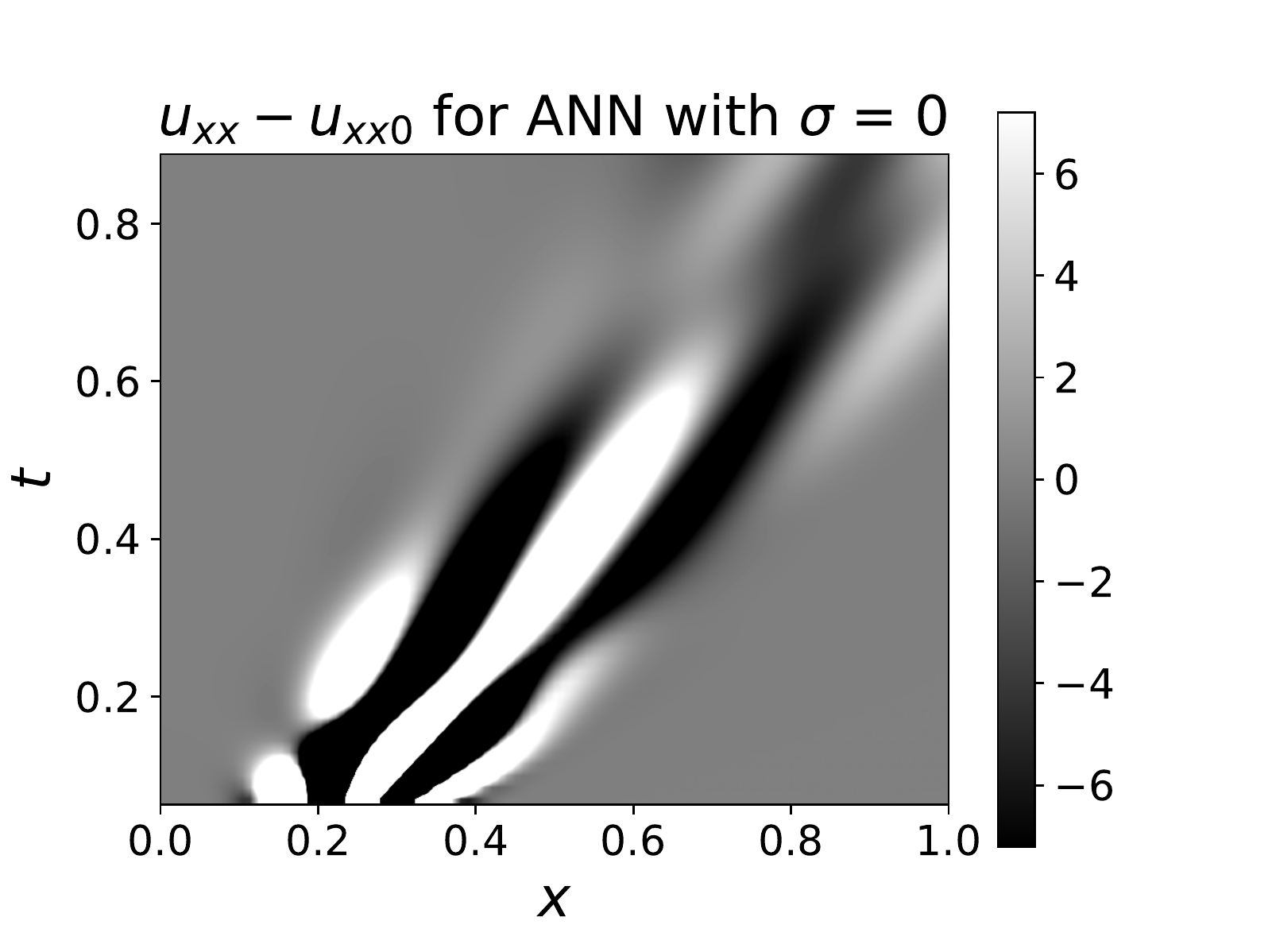}\includegraphics[width=0.4\textwidth]{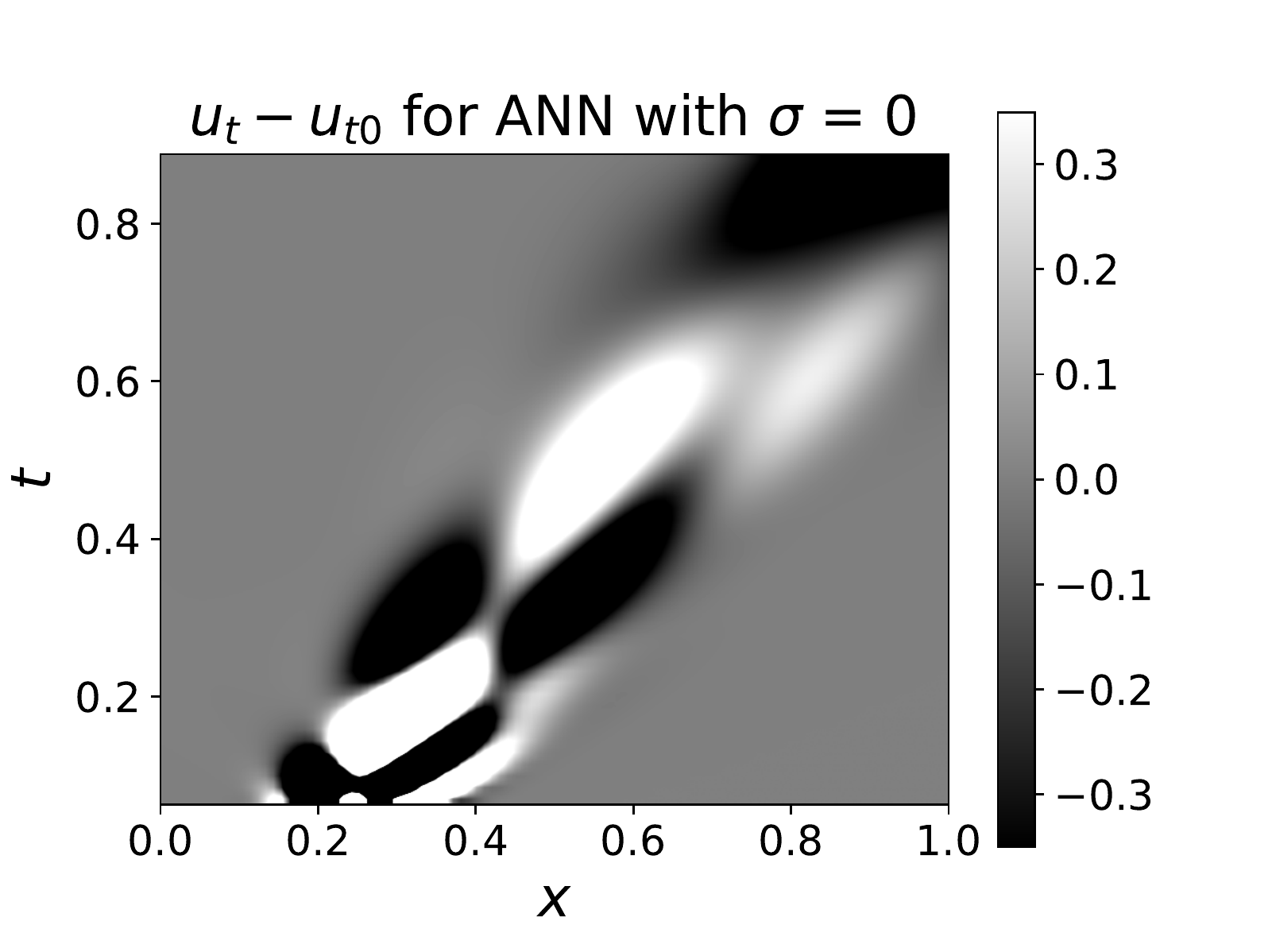}
    \caption{Contours of the residual values between the ANN model ($u, u_x,$, etc.) and analytical solutions ($u_0, u_x0,$ etc.) for the diffusion-advection data with $\sigma=0.$ Top left: Residuals for $u$, top right: residuals for $u_{x}$, bottom left: residuals for $u_{xx},$ bottom right: residuals for $u_{t}$. }
    \label{ANN_residual}
\end{figure}

We found that, in general, the PDE-FIND method learns the wrong equation, even when no noise is added to the data (Supplementary section S2, Figures S4-S6).  Each denoising method resulted in accurate estimates for $u(x,t)$ and its partial derivatives in this case, however. For example, the residuals between the ANN model and the analytical values for $u,u_t,u_x, \text{ and } u_{xx} $ were small  when $\sigma = 0$ (Figure \ref{ANN_residual}). We observed that, while small, the ANN residuals include systematic biases comprised of regions of over- and under-prediction. For example, all points near $(x,t) = (0.6,0.4)$ for the ANN's calculation for $u_x$ appear to over-predict the true value for $u_x$ in this region (Figure \ref{ANN_residual}). This contradicts the assumption of independence in $\{\epsilon_{i,j}\}_{i=1,j=1}^{M,N}$ for the  statistical model in Equation \eqref{stat_model}. As we will now demonstrate, these small, systematic error terms from the ANN cause PDE-FIND to learn the incorrect equation.

We illustrate here that PDE-FIND learns the incorrect equation when training data for $u_t$ is comprised of $u_t$ at all spatial points for the first half of the given time points and the validation data is comprised of all spatial points for the second half of all time points. Recall that in our actual implementation discussed below, we randomly split the training and validation data in $5\times5$ bins of adjacent spatiotemportal points. Using denoised values for $u(x,t)$ and its partial derivatives from the ANN in the case where $\sigma=0$ in the data, our training-validation procedure without pruning learns an equation of the form

\begin{equation}
u_{t}=a +bu_{x}+cu_{xx}+ du^2 + eu + fu^2u_x + gu^2u_{xx},\ a,...,g\in\R.\label{PDEFind_result}
\end{equation}
Similarly, the learned equations using finite difference and spline computations are

\begin{equation}
u_{t}=au_{x}+bu_{xx}+cu^{2}u_{x},\ a,b,c\in\R\label{PDEFind_result_fd}
\end{equation}
and 
\begin{equation}
u_{t}=au_{x}+bu_{xx}+cu^2u_x + du^2u_{xx} + eu_x^2,\ a,...,e\in\R,\label{PDEFind_result_spline}
\end{equation}
respectively.

\begin{figure}[h!]
    \centering
    \includegraphics[width=0.4\textwidth]{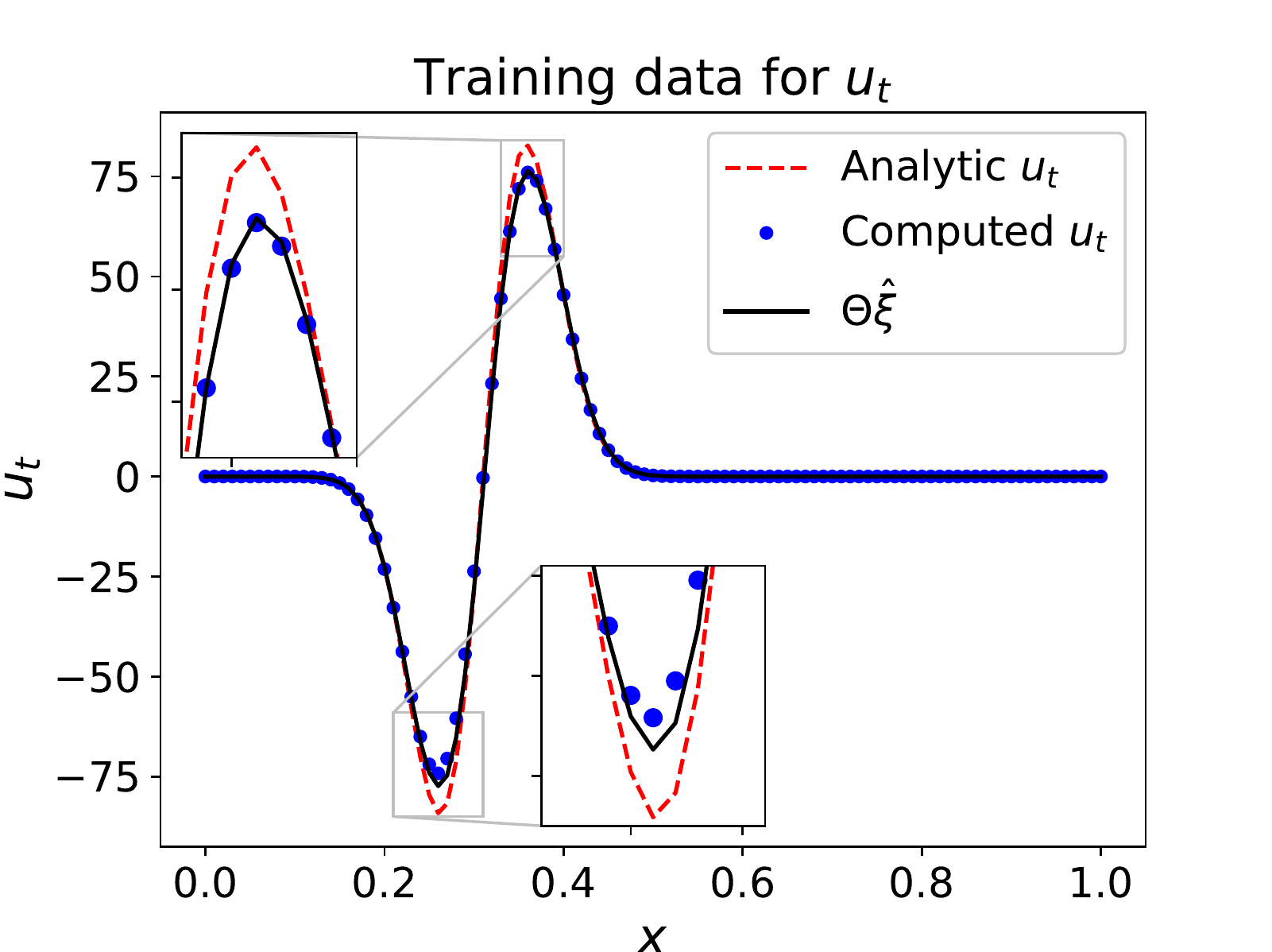}\includegraphics[width=0.4\textwidth]{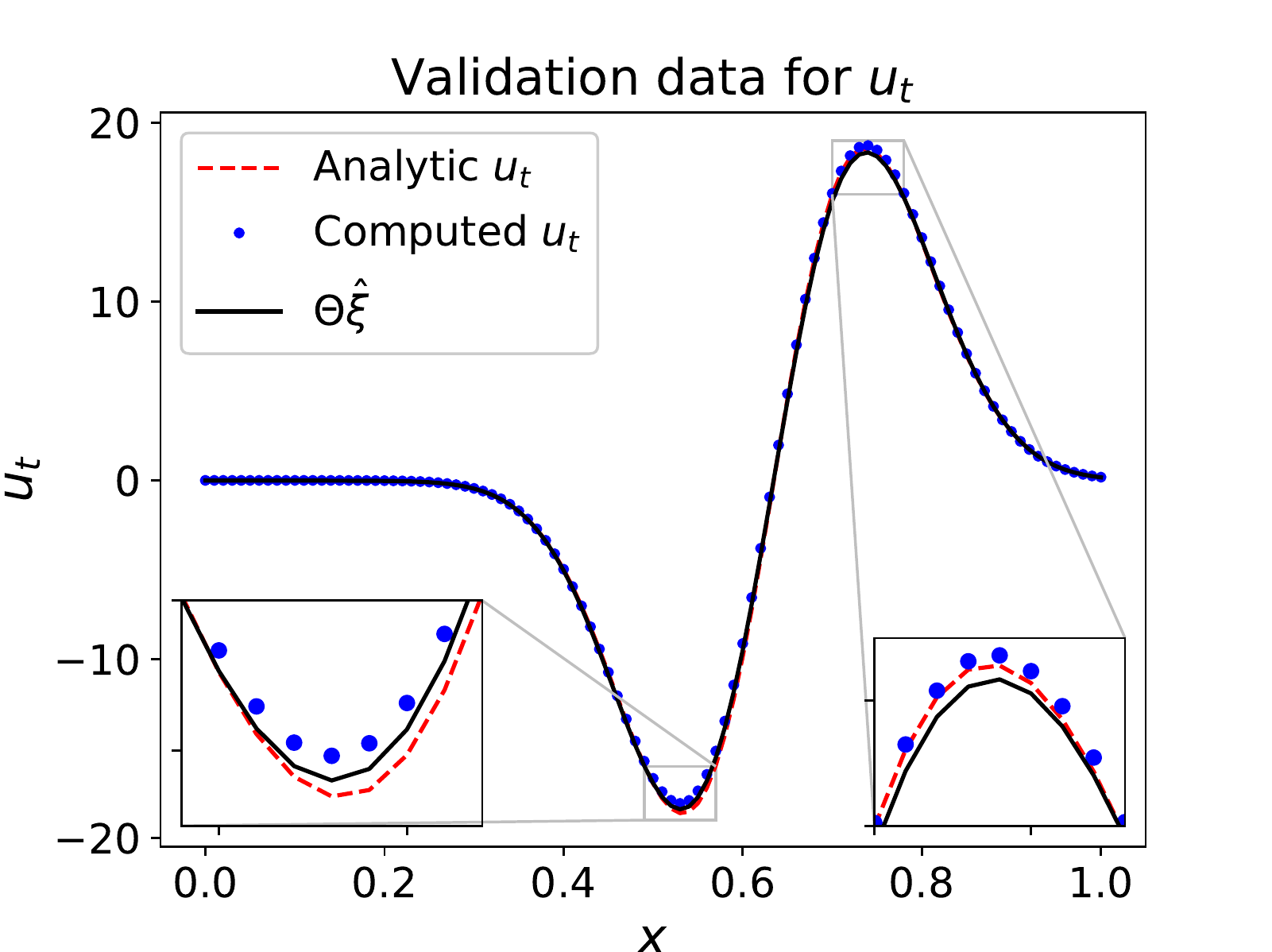}\caption{Results from the training (left) and validation (right) procedures
    in PDE-Find. The Red dashed line denotes the analytical value of $u_{t}$,
    the blue dots denote the computed $u_{t}$ values from the ANN, and the black lines denotes the equation for
    $u_{t}$ that has been computed with PDE-FIND. \label{train_val}}
\end{figure}

Each of these equations are incorrect and have extra terms on the right hand side of the learned PDE for the diffusion-advection equation. In Figure \ref{train_val}, we depict illustrative portions of
the training and validation sets comparing %
the analytical values of $u_{t}$ against the computed values of $u_t$ and PDE-FIND's selected equation using ANN approximations.  We found that PDE-FIND selects Equation \eqref{PDEFind_result} in place of the true diffusion-advection equation because
it recovers the ANN's incorrect computations of $u_t$ in both the training and validation data. In doing so, PDE-FIND fits the erroneous $u_{t}$ computations from the ANN approximation by including  extra terms in the learned PDE.

\subsection{PDE-FIND with pruning for the diffusion-advection equation} \label{PDE_find_DA}
\begin{figure}[h]
    \centering
    \includegraphics[width=.95\textwidth]{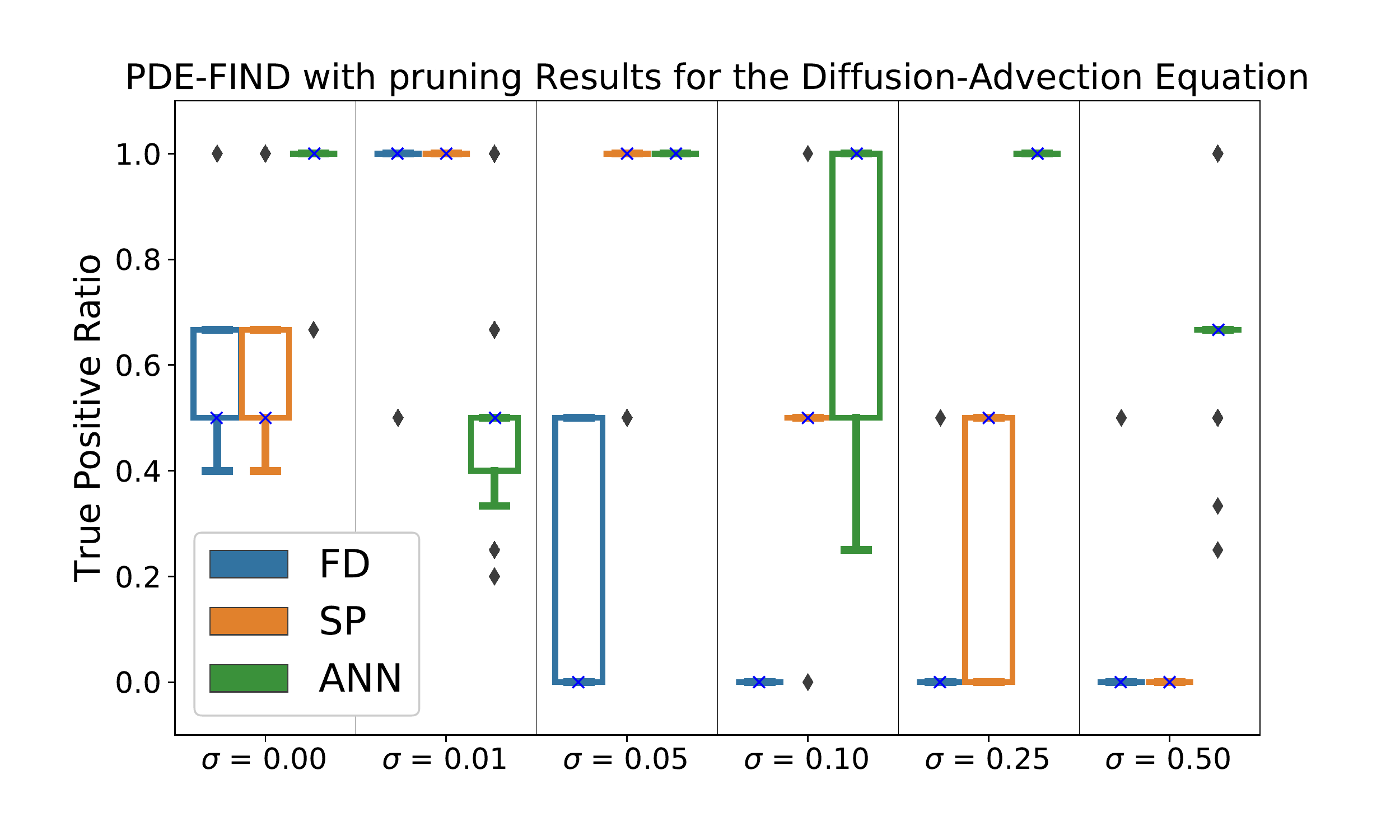}
    \caption{TPR values for the diffusion-advection equation. We calculated the TPR, see Equation \eqref{TPR}, for 1,000 different training-validation splits. These plots demonstrate the range of TPR values for each case. In each plot, the lower line in the colored box portion provides the 25\% quartile of the data and the upper line denotes the 75\% quartile. The ``x'' on each box plot denotes the median TPR value for that scenario. The length of the upper and lower whiskers are 1.5 times the interquartile range of the distribution, and diamonds denote outlier points. Any plot depicted as a solid horizontal line (e.g., the neural net computations for $\sigma=0$) denotes that that this value is the majority of the range of the distribution.}
    \label{box_whisker_DA}
\end{figure}

We tested whether implementing an additional pruning step with PDE-FIND could remove the extra terms
resulting from the biases discussed in Section \ref{PDE-FIND_pruning}. For the diffusion-advection equation, we found that for all values of $\sigma$ except $\sigma=0.01$, PDE-FIND with pruning achieves the highest median TPR when using ANN approximations (Figure \ref{box_whisker_DA}). The ANN's median value is TPR = 1 (meaning that over half of the simulations yielded the correct equation form) for $\sigma=0,0.05,0.10,\text{ and }0.25$. The ANN resulted in a median TPR = 0.667 at $\sigma=0.50$. In contrast, the spline method only achieved of median TPR = 1 at the lower noise levels $\sigma=0.01 \text{ and } 0.05$. For $\sigma\ge0.10$, the medians for the spline method were all TPR$\leq$0.5. The finite difference method resulted in a median TPR = 1 at $\sigma=0.01$, but the median TPR = 0 for larger values of $\sigma$.

Supplementary Table S1 shows the most commonly learned PDEs for each denoising method at each noise level. We found that the ANN method, used in conjunction with PDE-FIND with pruning, resulted in the correct PDE for $\sigma=0,.05,0.10,0.25$. The ANN specifies the incorrect equation for $\sigma=0.01$ and $\sigma=0.50$. However in both of these cases, the extra terms have small parameter values (e.g. 0.001) that a scientist with an understanding of the system under consideration may manually neglect. On the other hand, PDE-FIND cannot discover the correct equation with finite difference or spline computations for $\sigma\ge0.10$. These results suggest that the ANN method enables PDE-FIND with pruning to learn the diffusion-advection equation accurately at biologically realistic noise levels, e.g, $\sigma=.05,0.10,0.25$.

\subsection{PDE-FIND with pruning for the Fisher-KPP Equation} \label{PDE_find_FKPP}

\begin{figure}[h]
    \centering
    \includegraphics[width=.95\textwidth]{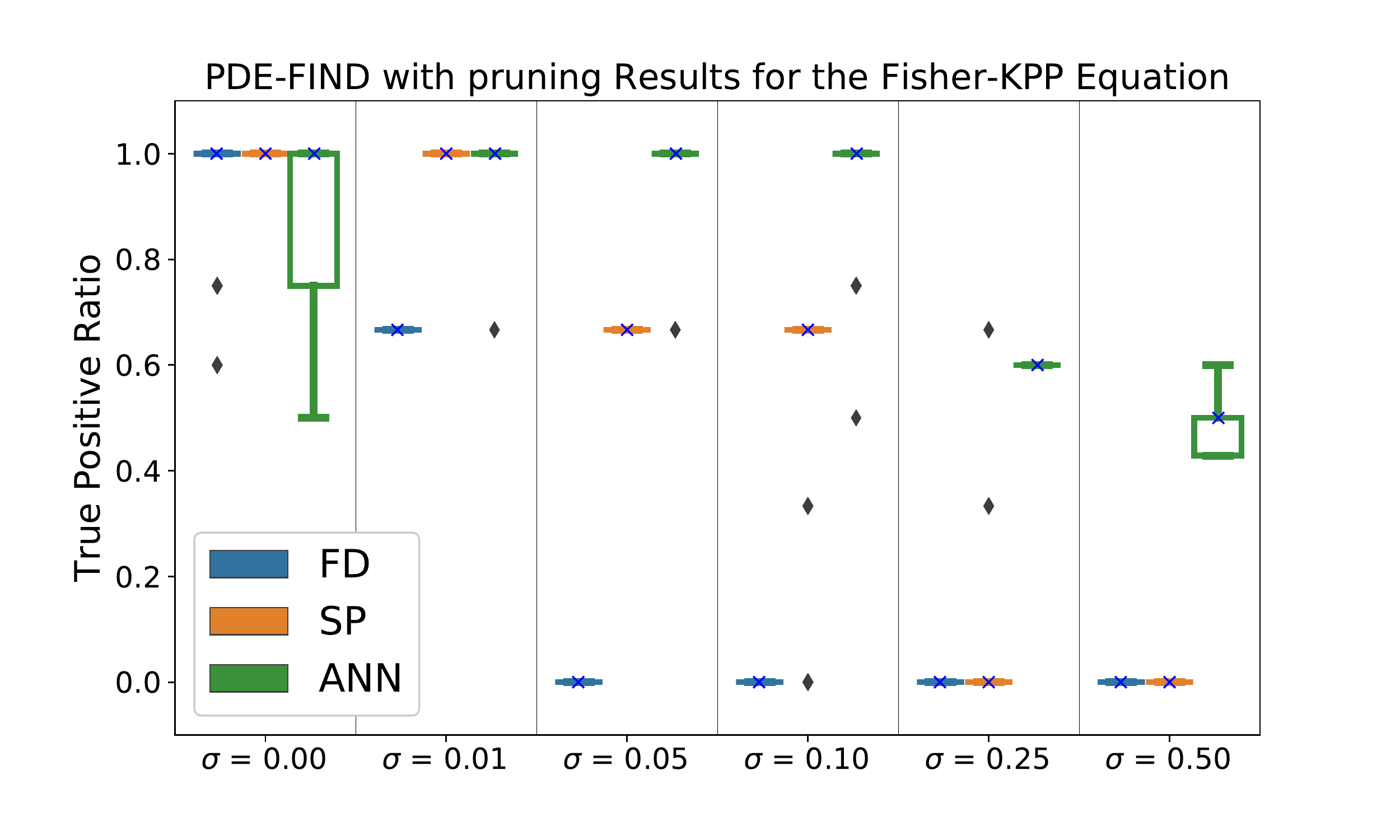}
    \caption{TPR values for the Fisher-KPP Equation. We calculated the TPR, see Equation \eqref{TPR}, for 1,000 different training-validation splits. These plots demonstrate the range of TPR values for each case. In each plot, the lower line in the colored box portion provides the 25\% quartile of the data and the upper line denotes the 75\% quartile. The ``x'' on each box plot denotes the median TPR value for that scenario. The length of the upper and lower whiskers are 1.5 times the interquartile range, and diamonds denote outlier points. Any plot depicted as a solid horizontal line (e.g., the finite difference computations for $\sigma=0$) denote that this value is the majority of the distribution. }
    \label{box_whisker_fisher}
\end{figure}

We tested the PDE-FIND with pruning method in conjuction with several denoising strategies using data from the Fisher-KPP Equation. We found that the ANN method had a median TPR = 1.0 (meaning that the correct equation is specified for at least half of the training-validation data splits) for $\sigma=0,0.01,0.05,\text{ and }0.10$ (Figure \ref{box_whisker_fisher}). In contrast, the finite difference calculations only had a median TPR = 1 at $\sigma = 0 $, and the spline method only had median TPR = 1 at $\sigma=0,0.01$. The accuracy in using PDE-FIND with the spline and finite difference methods quickly deteriorates for high noise levels. The finite difference method resulted in a median TPR = 0 for $\sigma\ge0.05$, and the spline method resulted in a median TPR = 0.667 at $\sigma=0.05,0.10$ and TPR = 0 at $\sigma=0.25,0.50$. The ANN had a median TPR = 0.6 and 0.5 for $\sigma=0.25\text{ and } 0.5,$ respectively.

Supplementary Table S2 shows the most commonly chosen PDEs resulting from the PDE-FIND with pruning method. We found that PDE-FIND with pruning is able to discover the correct equation form for $\sigma=0,0.01, 0.05,\text{ and }0.10$ when using the ANN approximations. While PDE-FIND is unable to specify the correct equations with ANN data for $\sigma=0.25$ and 0.50, all of the terms in the Fisher-KPP equation were included in the learned PDEs. In contrast, using the spline method for denoising resulted in only learning the correct PDE for $\sigma\leq0.01$. For larger values of $\sigma$, the spline method resulted in large errors in the derivative approximations and did not yield any terms on the right hand side of the learned PDE. Similarly, the finite difference method resulted in only learning the true equation form for $\sigma=0$. These results suggest that only the ANN method enables PDE-FIND with pruning to learn the Fisher-KPP equation for reasonably high noise levels of $\sigma=0.05,0.10$.

\subsection{PDE-FIND with pruning for the nonlinear Fisher-KPP Equation} \label{PDE_find_nonlinear_FKPP}

\begin{figure}[h]
    \centering
    \includegraphics[width=0.95\textwidth]{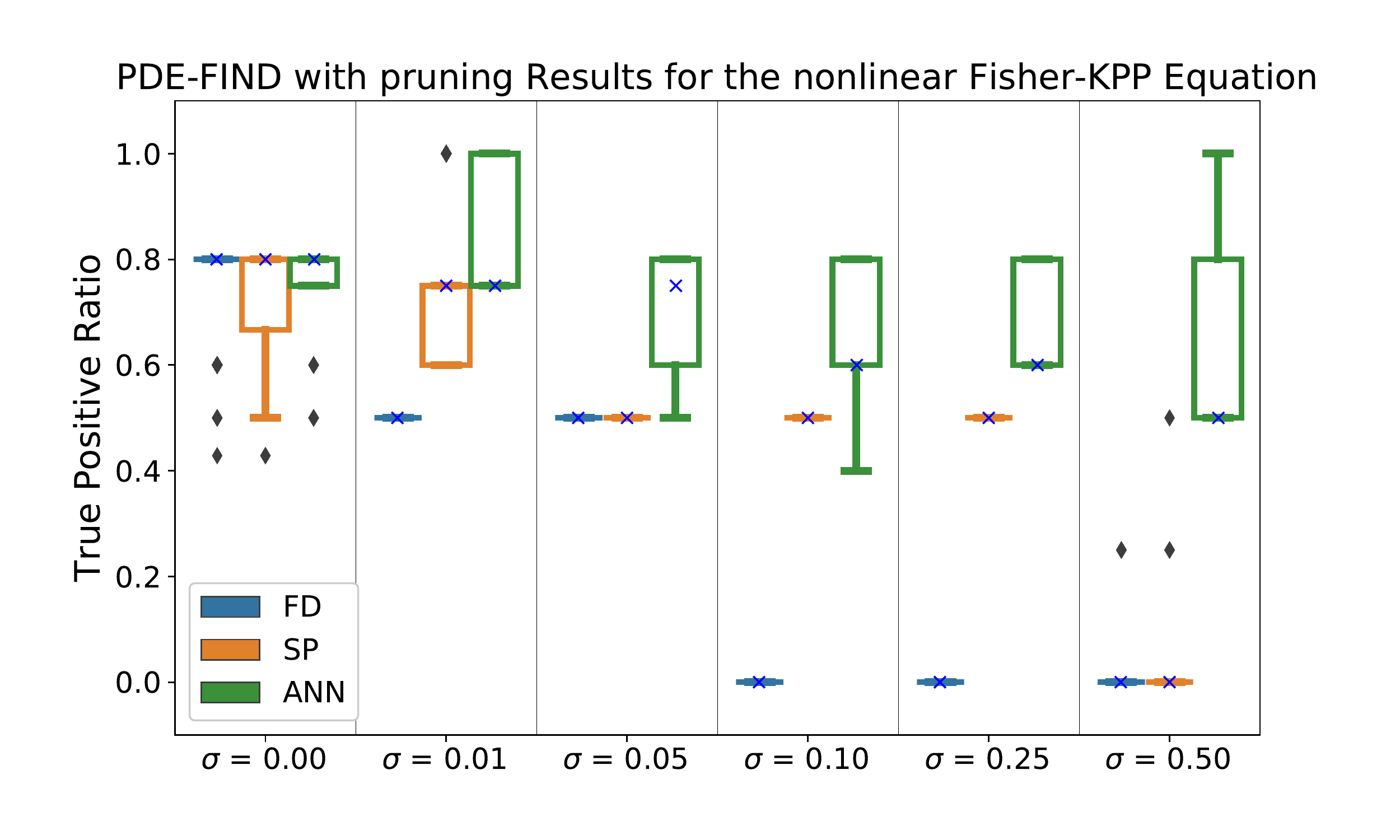}
    \caption{TPR values for the nonlinear Fisher-KPP Equation. We calculated the TPR, see Equation \eqref{TPR}, for 1,000 different training-validation splits. These plots demonstrate the range of TPR values for each case. In each plot, the lower line in the colored box portion provides the 25\% quartile of the data and the upper line denotes the 75\% quartile. The ``x'' on each box plot denotes the median TPR value for that scenario. The length of the upper and lower whiskers are 1.5 times the interquartile range, and diamonds denote outlier points. Any plot depicted as a lone solid horizontal line (e.g., the finite difference computations for $\sigma=0$) denotes that that this value is the whole range of the data. }
    \label{box_whisker_fisher_nonlin}
\end{figure}

We found that the PDE-FIND with pruning method struggles with all three denoising strategies to recover the correct PDE from data that has been generated by the nonlinear Fisher-KPP Equation. PDE-FIND could not achieve a median TPR = 1 for any of these methods, meaning that the correct equation was never specified for over half of the training-validation data splits (Figure \ref{box_whisker_fisher_nonlin}). All methods have median TPR = 0.8 at $\sigma=0$. When using ANN approximations, the PDE-FIND with pruning method has median TPR = 0.8 at $\sigma=0.01\text{ and }0.05$, TPR = 0.6 at $\sigma=0.10\text{ and }0.25$, and TPR = 0.5 at $\sigma=0.50$. When using spline computations, PDE-FIND with pruning has median TPR = 0.8 at $\sigma=0.01$, TPR = 0.5 at $\sigma=0.05,0.10,\text{ and }0.25$, and TPR = 0 at $\sigma=0.50$. When using finite difference computations, PDE-FIND with pruning has median = 0.50 at $\sigma=0.01\text{ and }0.05$ and TPR = 0 for $\sigma\ge0.10$.  

While all of the denoising strategies lead to incorrect equations, the ANN strategy recovers the most relevant terms in its final equations. Supplementary Table S3 shows the most-commonly chosen PDEs resulting from the PDE-FIND with pruning method. We found that all three methods predict the true equation with an extra Fickian diffusion term, $u_{xx}$, for $\sigma = 0$. When using the ANN approximations, the PDE-FIND with pruning algorithm recovers the true equation with added Fickian diffusion at $\sigma=0.05$, and it recovers three of the correct terms but excludes $uu_{xx}$ at $\sigma=0.01$. For larger values of $\sigma$ with the ANN approximations, PDE-FIND with pruning recovers three correct terms, excludes the $uu_{xx}$ term, and includes an extra Fickian diffusion term (as well as an additional constant term at $\sigma=0.50$). When using spline computations, the PDE-FIND with pruning algorithm  recovers three correct terms, excludes the $uu_{xx}$ term, and adds an extra Fickian diffusion term at $\sigma=0.01$. For $\sigma=0.05-0.25$, the final equation recovers two correct terms but excludes the $uu_{xx}$ and $u_x^2$ terms. At $\sigma=0.50,$ all terms are deleted when using spline computations. When using finite difference computations, the PDE-FIND with pruning algorithm correctly recovers two terms but excludes the $uu_{xx}$ and $u_x^2$ terms at $\sigma=0.01\text{ and }0.05$. For larger values of $\sigma$, no correct terms are included in the final equation form. 

We investigated if the recovered terms from the PDE-FIND with pruning algorithm using ANN approximations can be used as the specified mathematical model in an inverse problem methodology (cf., \cite{banks_mathematical_2009}) to recover the final parameter estimate values from the nonlinear Fisher-KPP Equation. If we take the union of all terms that are included in the final equations in Supplementary Table S3 for the ANN method using noisy data ($\sigma>0$), then we have an equation of the form
\begin{equation}
    u_t = au_{xx} + buu_{xx} + cu_x^2 + du + eu^2 + f, a,\hdots,f\in\R. \label{nonlin_FKPP_learned_equation}
\end{equation}
We estimated the parameters $a,...,f$ in Equation \eqref{nonlin_FKPP_learned_equation} for each value of $\sigma$ by simulating the solution to this PDE using the method of lines and minimizing Equation \eqref{cost} using the Nelder-Mead algorithm. We input the equations from Supplementary Table S3 for the ANN method as the initial guess for each data set. We find that performing this inverse problem leads to accurate parameter estimates for the true terms in the nonlinear Fisher-KPP Equation and small coefficient values for the incorrect terms ($u_{xx}$ and $1$)  for $\sigma=0,0.01,0.05,\text{ and }0.25$  (Table \ref{nonlin_FKPP_IP}). At $\sigma=0.10$ and $0.50$, this inverse problem methodology leads to small coefficient estimates for $uu_{xx}$ in addition to $u_{xx}$ and $1$. Note that this same process would not lead to ultimately recovering the true equation and parameter estimates from the spline or finite difference approximations because their final equations never included the correct $uu_{xx}$ term in the final equation for noisy data ($\sigma>0$).

\begin{table}[h]
    \centering
    \caption{Inferred parameters for the nonlinear Fisher-KPP Equation data when performing an inverse problem on Equation \eqref{nonlin_FKPP_learned_equation}.}
    \begin{tabular}{|c|l|}
    \hline 
     &  \textbf{True Equation}\tabularnewline
    \hline 
     & $u_{t}=.02uu_{xx}+.02u_{x}^{2}+10u-10u^{2}$\tabularnewline
    \hline 
    $\sigma$ & \textbf{Revised Equation}\tabularnewline
    \hline 
    0 & $u_{t}=2.3\times10^{-12}u_{xx}+.017uu_{xx}+.021u_{x}^{2}+10.0u-10.0u^{2} -1.60\times10^{-8}$\tabularnewline
    \hline 
    0.01 & $u_{t}=1.4\times10^{-5}u_{xx}+.021uu_{xx}+.019u_{x}^{2}+10.0u-10.0u^{2}-1.38\times10^{-8}$\tabularnewline
    \hline 
    0.05 & $u_{t}=2.5\times10^{-4}u_{xx}+.0034uu_{xx}+.029u_{x}^{2}+9.9u-10.0u^{2} + -2.86\times10^{-8}$\tabularnewline
    \hline 
    0.10 & $u_{t}=6.1\times10^{-4}u_{xx}+8.42\times10^{-4}uu_{xx}+.023u_{x}^{2}+9.38u-9.52u^{2}-2.44\times10^{-4}$\tabularnewline
    \hline 
    0.25 & $u_{t}=7.2\times10^{-4}u_{xx}+.015uu_{xx}+.024u_{x}^{2}+9.8u-9.14u^{2}+3.53\times10^{-7}$\tabularnewline
    \hline 
    0.50 & $u_{t}=2.1\times10^{-3}u_{xx}-3.72\times10^{-3}uu_{xx}+.032u_{x}^{2}+9.7u-7.5u^{2}+1.38\times10^{-6}$\tabularnewline
    \hline 
    \end{tabular}
    \label{nonlin_FKPP_IP}

\end{table}

\section{Conclusions and Future Work} \label{Conclusions}

The novel use of the ANN method presented here is a significant step toward making PDE learning more achievable in realistic scenarios with noisy biological data. Because an ANN is a fully differentiable function, it can be used to approximate derivative computations to build the library of terms needed for PDE learning. The current practice to build a library of terms for learning PDEs when noise is present in the observed data $u(x,t)$ is to use finite difference or spline approximations for small amounts of noise \cite{rudy2017data,zhang2018robust}. Our findings suggest that these two methods are highly sensitive to the amount of noise in the data in the range of less than $5\%$ noise. Moreover, these methods can not incorporate heteroscedastic statistical models such as Equation \eqref{stat_model} to handle non-constant error noise, which is a typical phenomenon encountered in biological data. We showed that the ANN method outperforms spline and finite difference approximations of $u(x,t)$ and its partial derivatives when significant levels of non-constant error noise are present in the data $u(x,t)$. 

In general, we found that when using the ANN method, the true underlying equation is more accurately recovered than when using finite differences or splines at all noise levels considered in this work (up to $\sigma=0.50$). The disparity between the accuracy of the ANN and the other two methods sharply increased with the noise levels in the data. For the diffusion-advection equation, the ANN allows PDE-FIND to learn the correct equation for up to 25\% noise levels, while splines and finite difference calculations fail for noise levels over 10\%. The ANN data allows PDE-FIND to learn the Fisher-KPP Equation for up to 10\% noise levels whereas the other two methods fail at 5\%. All methods struggle with the nonlinear Fisher-KPP Equation, but the ANN data leads to an equation that can then be used with an inverse problem methodology \cite{banks_mathematical_2009} to infer which terms are meaningful.

There are several reasonable explanations for the difficulties in learning the nonlinear Fisher-KPP Equation. Three of the terms in Equation \eqref{nonlinear_fisher} are the product of two terms including $u$ or its derivatives ($u^2$,$u_x^2$, and $uu_{xx}$). In practice, these terms may be inaccurate approximations from noisy data (as demonstrated in Table \ref{fisher_nonlin_rmse}). Multiplying two inaccurate terms may lead to an even larger amount of uncertainty associated with these estimates. We postulate that the high level of uncertainty in these type of terms resulting from the product of inaccurate estimates likely increases the difficulty of learning to include them in the process of PDE learning. %
Furthermore, it must be noted that the data for this equation was generated numerically the by finite difference method. Though some analytical solutions to the Fisher-KPP equations are known, they come either in series form, which introduces error by series truncation, or in traveling wave form, which would be indistinguishable from the advection equation to PDE-Find. The finite difference methods used to approximate the spatial derivatives in the nonlinear Fisher-KPP equation introduce second and fourth-order truncation error terms which lead to numerical dispersion and diffusion effects %
that may account for the recovery of some unexpected terms. 

We found that the use of pruning following our implementation of the PDE-FIND algorithm increased our ability to recover the correct equation in terms of the TPR. It may be argued that these additional terms learned from the PDE-FIND implementation without pruning would be removed from the final equation if more regularization (i.e., a larger value of $k$ in implementation of the Greedy algorithm) were used. However, when faced with the issue of learning the governing equation from actual data in practice, one will not have the ability to know when the specified equation is correct or not. We thus need to identify the correct hyperparameters without any \emph{a priori} knowledge. We observed that the systematic biases in our ANN (depicted in Figure \ref{ANN_residual}) make it difficult to choose a hyperparameter value that leads to the correct equation because these biases are present in the training and validation data. The pruning algorithm is a way to correct for when the incorrect hyperparameter has been chosen by ensuring that all terms in the learned PDE are sufficiently sensitive to constitute a strong signal in the data, instead of resulting from a bias in our approximation methods.

Our use of pruning to remove terms from learned PDEs could be improved in future work. While effective, our implementation is somewhat crude, in which we pre-specify a threshold level to prune parameters based on out-of-sample MSE values on the validation data set. Previous studies have discussed F-statistics as one way to infer the increase in variance that pruning a variable will lead to, but there are still many different interpretations of these results which makes a definitive statistical pruning method challenging to ascertain  \cite{burgess_non-linear_1995}.

In future work, we may also consider extending the use of the ANN method to statistical models that a relevant to other types of biological noise, such as weighted least squares statistical models \cite{banks_mathematical_2009} or log-normally distributed errors \cite{bortz_model_2006}. We will also investigate using our methodology in the case where the statistical model is not known and needs to be estimated from data with a model-free numerical method \cite{banks_difference-based_2016}. %

\end{document}

\maketitle

\section{Comparing spline and bi-spline methods for denoising data}
\begin{figure}[h!]
    \centering
    \includegraphics[width=0.95\textwidth]{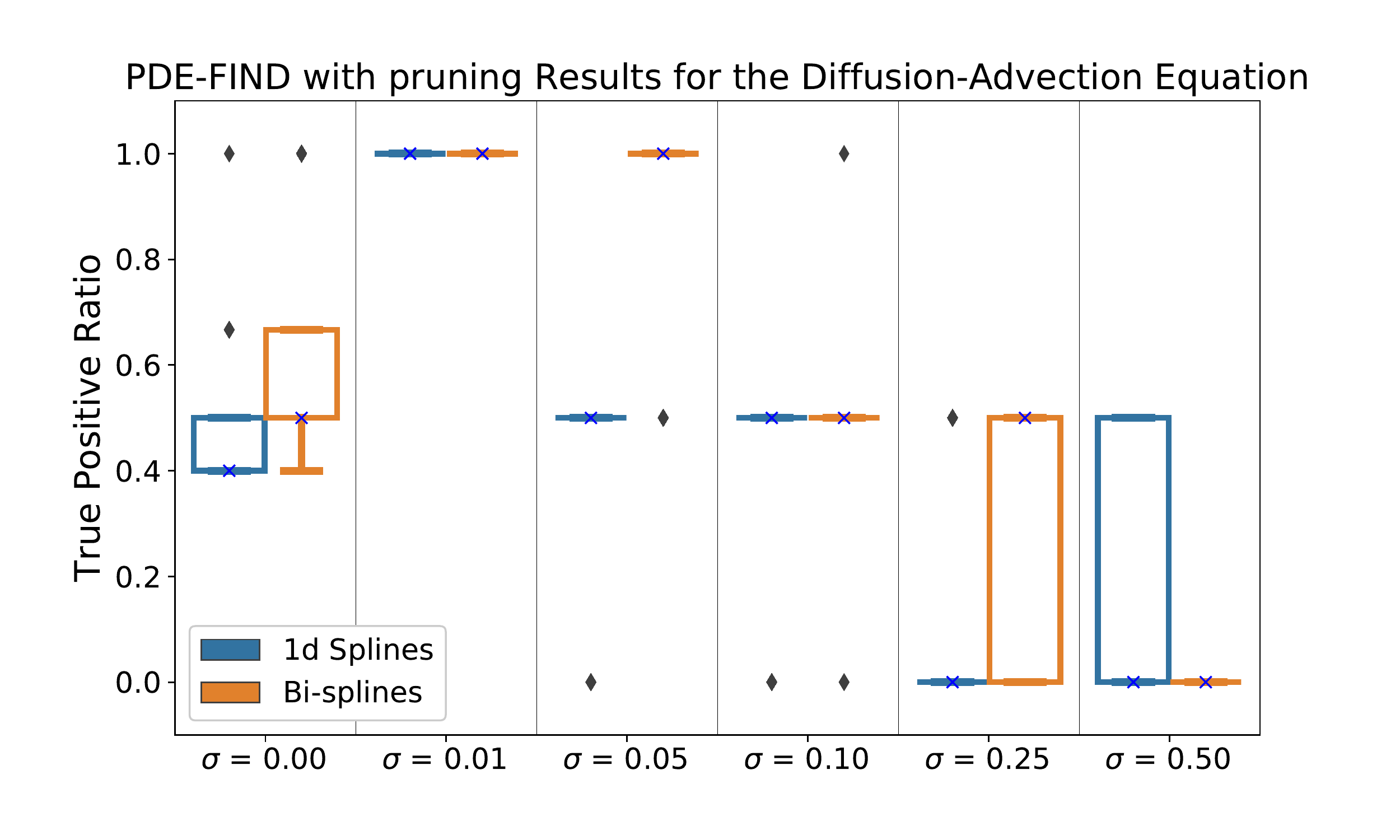}
    \caption{TPR values for the diffusion-advection equation when using 1-dimensional cubic splines versus cubic bi-splines for denoising data and approximating partial derivatives.}
    \label{box_whisker_diffadv_splines}
\end{figure}

\begin{figure}
    \centering
    \includegraphics[width=0.95\textwidth]{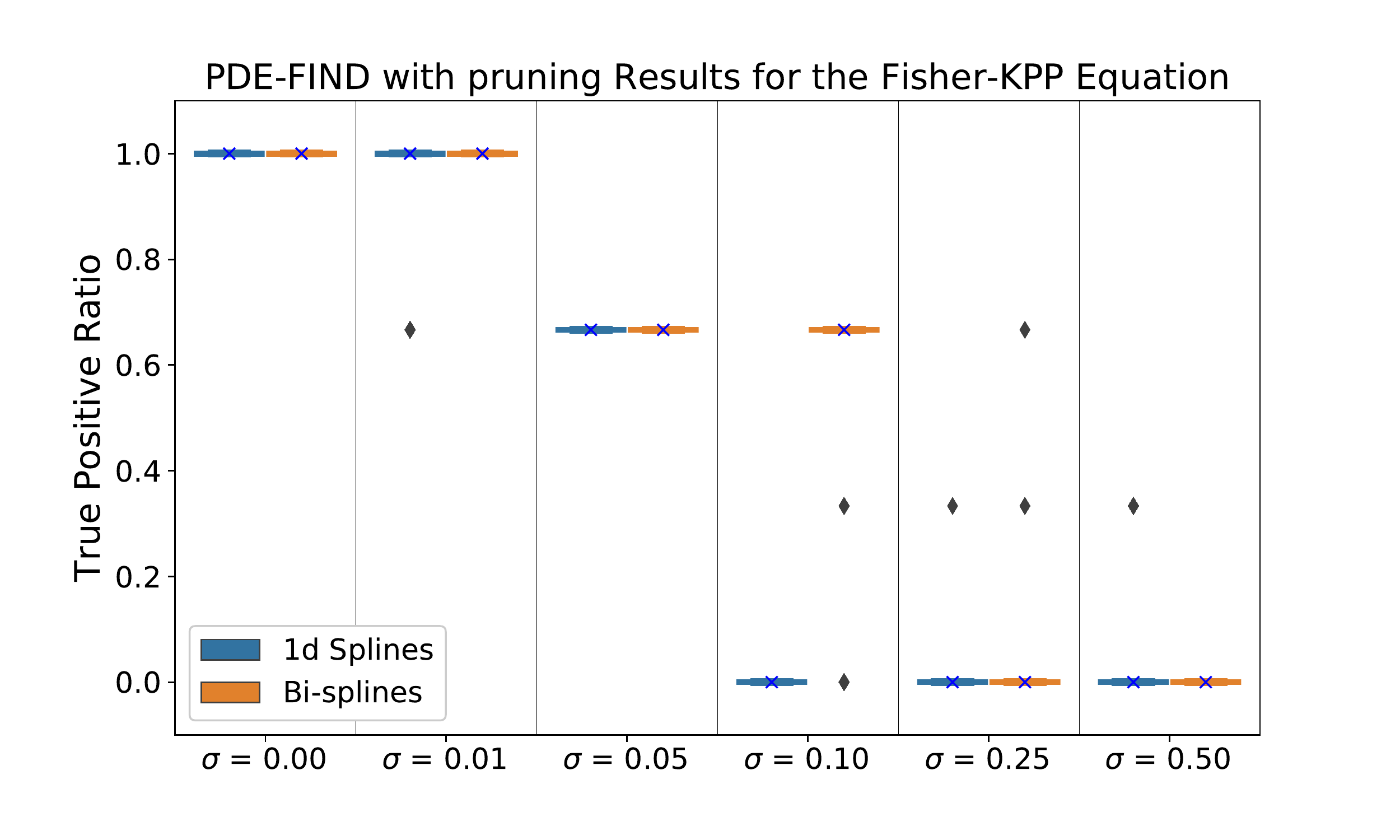}
    \caption{TPR values for the Fisher-KPP equation when using 1-dimensional cubic splines versus cubic bi-splines for denoising data and approximating partial derivatives. }
    \label{box_whisker_fisher_splines}
\end{figure}

\begin{figure}
    \centering
    \includegraphics[width=0.95\textwidth]{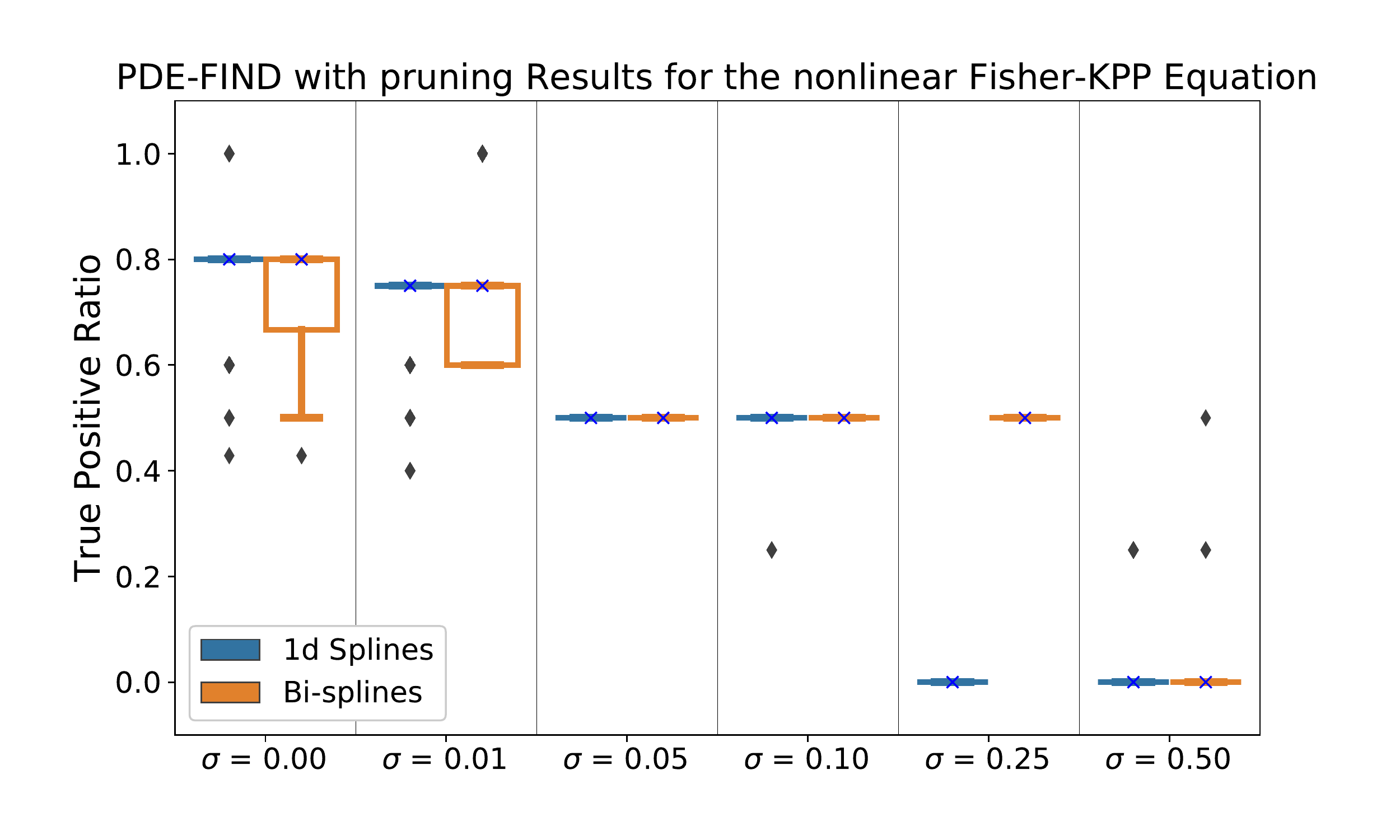}
    \caption{TPR values for the nonlinear Fisher-KPP equation when using 1-dimensional cubic splines versus cubic bi-splines for denoising data and approximating partial derivatives.  }
    \label{box_whisker_fisher_nonlin_splines}
\end{figure}

We compared the accuracy in learning the correct PDE when using 1-dimensional cubic splines versus cubic bi-splines for denoising data and approximating partial derivatives (Figures \ref{box_whisker_diffadv_splines},\ref{box_whisker_fisher_splines},\ref{box_whisker_fisher_nonlin_splines}). We found that PDE-FIND with pruning always has a higher TPR value when using bi-spline computations as compared to 1-dimensional splines.

\newpage

\section{PDE-FIND without pruning results}
\begin{figure}[h!]
    \centering
    \includegraphics[width=0.95\textwidth]{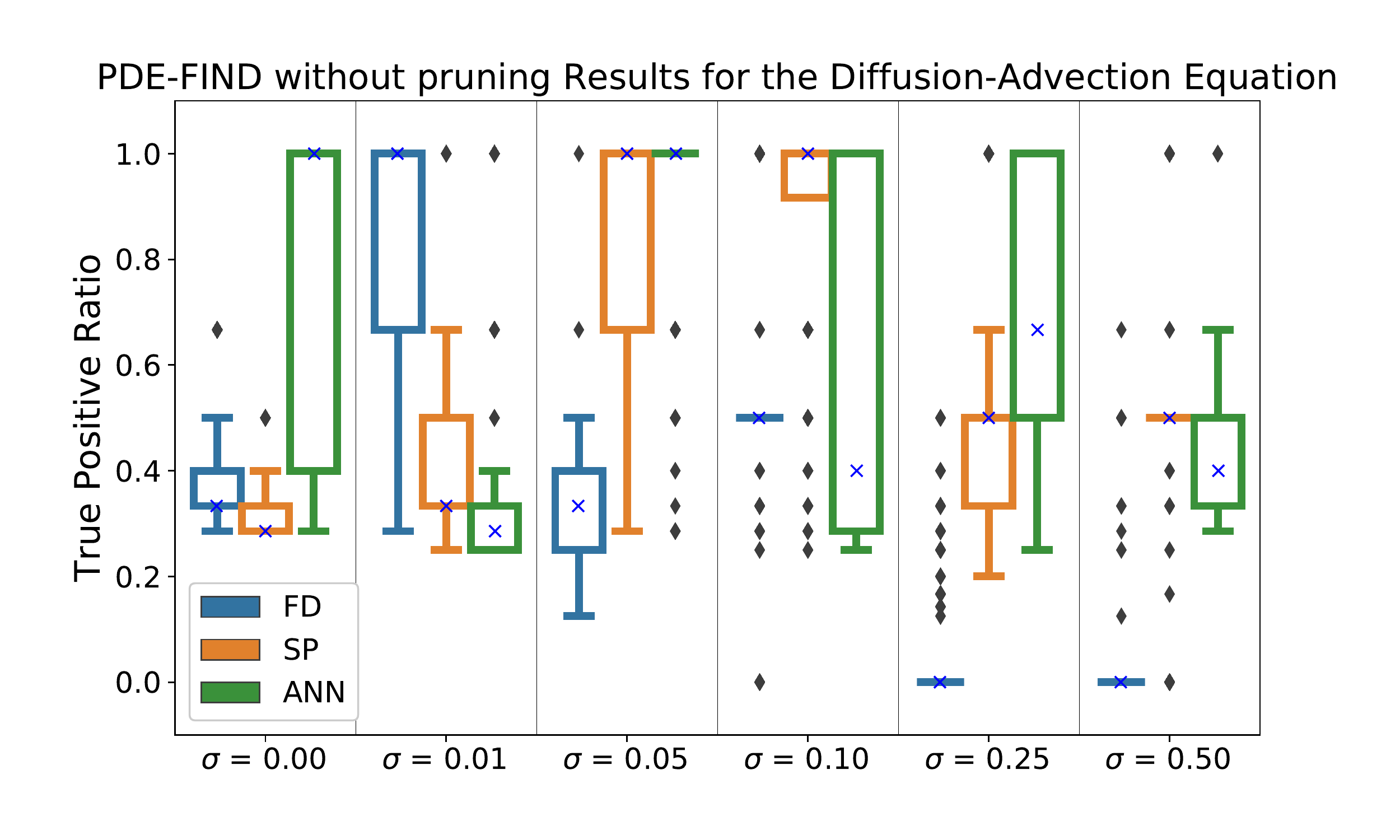}
    \caption{TPR values for the diffusion-advection equation.  }
    \label{box_whisker_diffadv_no_prune}
\end{figure}

\begin{figure}
    \centering
    \includegraphics[width=0.95\textwidth]{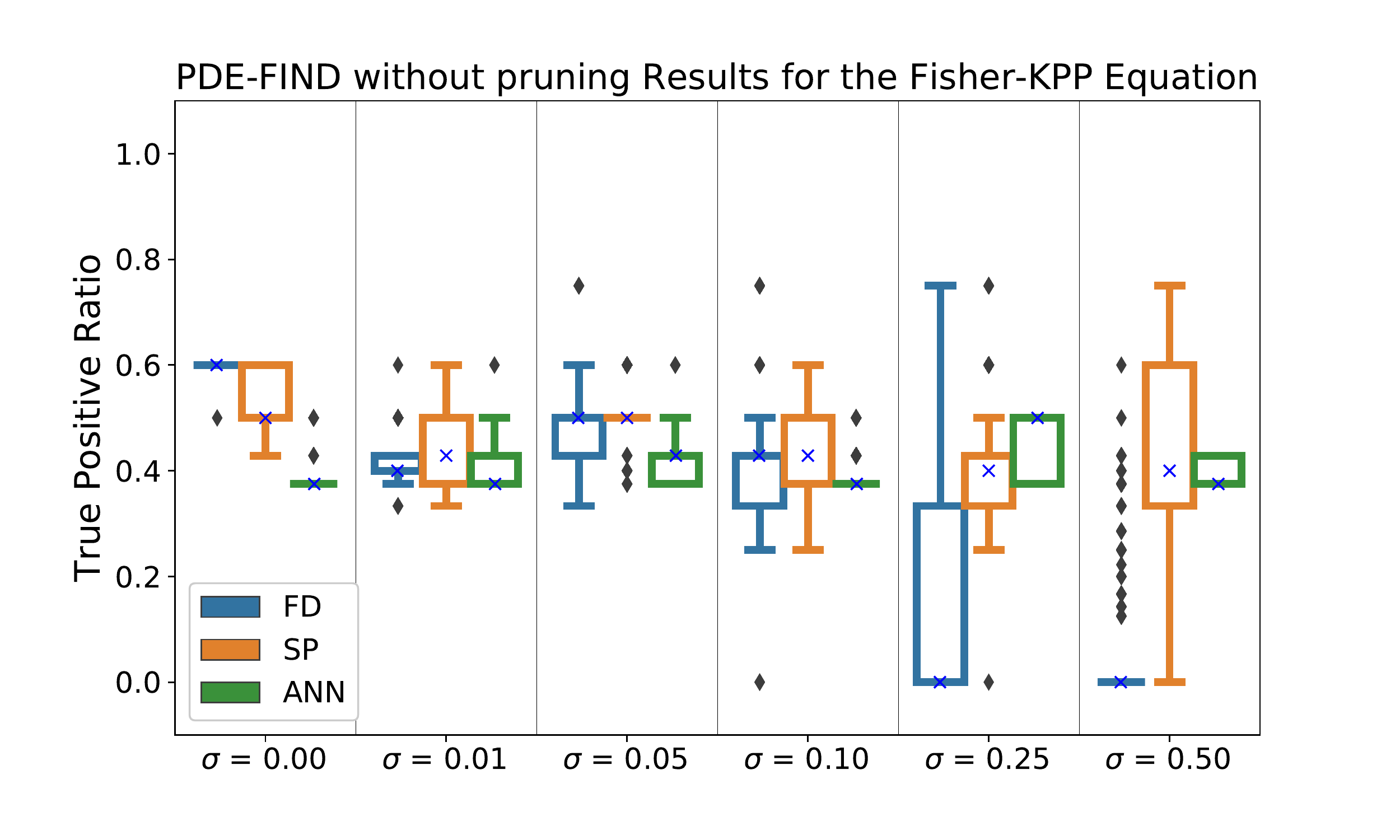}
    \caption{TPR values for the Fisher-KPP equation.}
    \label{box_whisker_fisher_no_prune}
\end{figure}

\begin{figure}
    \centering
    \includegraphics[width=0.95\textwidth]{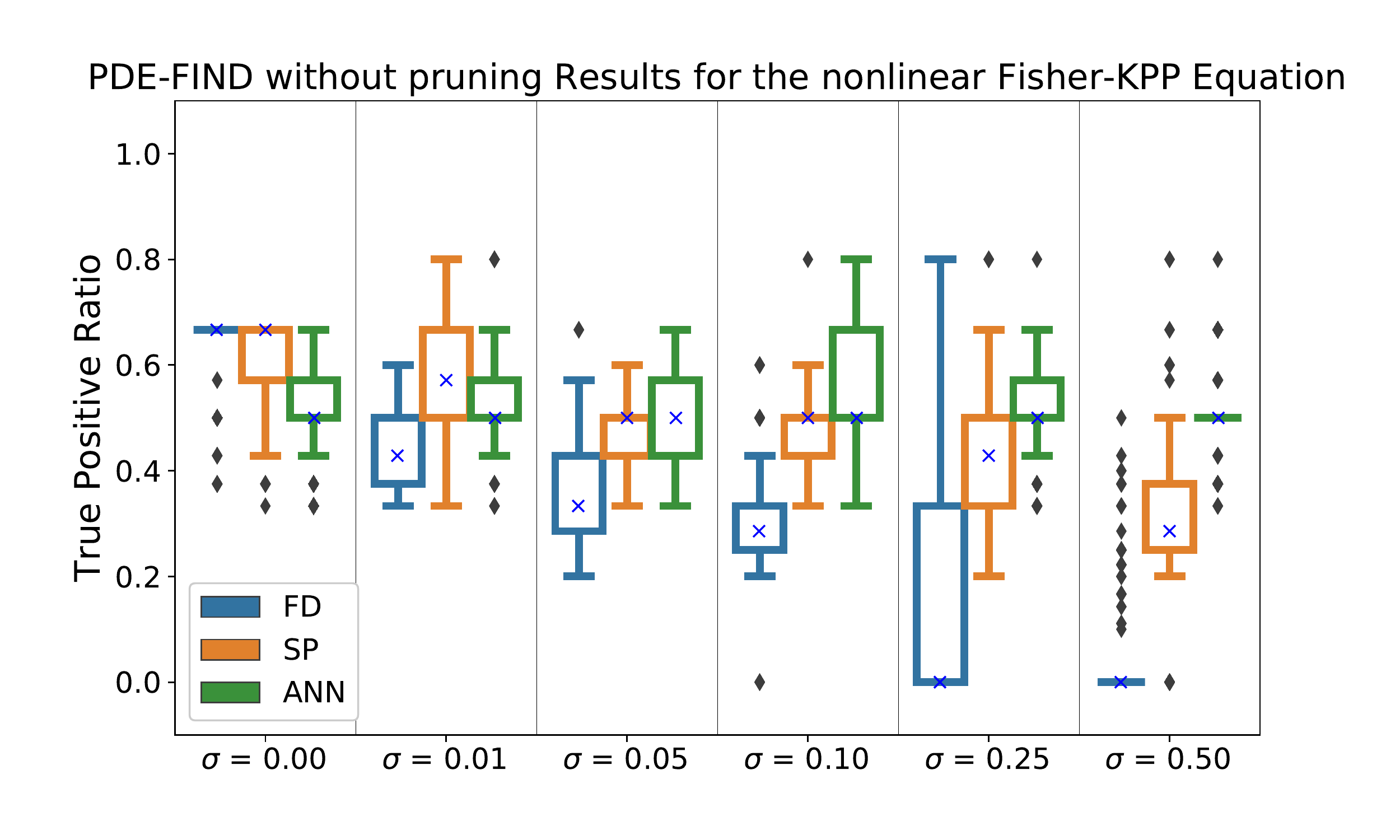}
    \caption{TPR values for the nonlinear Fisher-KPP equation.}
    \label{box_whisker_fisher_nonlin_no_prune}
\end{figure}

We found that using PDE-FIND without pruning results in learning the wrong equation when applied to data from biological transport models, even when no noise is added to the data. We evaluated accuracy, using the true positive ratio (TPR) as a metric, for the diffusion-advection (Figure \ref{box_whisker_diffadv_no_prune}), Fisher-KPP (Figure \ref{box_whisker_fisher_no_prune}), and nonlinear Fisher-KPP equations (Figure \ref{box_whisker_fisher_nonlin_no_prune}).

For the diffusion-advection equation, we found that the TPR value of the final learned equation when using ANN approximations is higher for all values of $\sigma$ when using pruning with PDE-FIND than without pruning (Figure \ref{box_whisker_diffadv_no_prune}). In general, for small values of $\sigma$, we observed that pruning enables PDE-FIND to better learn the true equation when using spline and finite difference computations, but it harms the ability to learn the true equation for larger values of $\sigma$. For example, the median TPR value increases after pruning when using finite difference approximations from TPR = 0.33 to 0.5 for $\sigma=0$. However, the TPR instead decreases from TPR = 0.33 to 0 at $\sigma=0.05$ and from TPR = 0.5 to 0 at $\sigma=0.10$. The median TPR value when using spline approximations increases from TPR = 0.3 to 0.5 at $\sigma=0$ and from TPR = 0.33 to 1 at $\sigma=0.01$. At $\sigma=0.10$, the median values decreased from TPR = 1.0 to 0.5. 

For the Fisher-KPP equation, the median TPR value for PDE-FIND with the ANN computations always increases after using pruning (Figure \ref{box_whisker_fisher_no_prune}). The median value for PDE-FIND with finite difference computations increases for $\sigma=0,0.01$, but decreases from TPR = 0.5 to 0 for $\sigma=0.05$ and from TPR = 0.45 to 0 for $\sigma=0.10$. The median TPR value for PDE-FIND with spline computations increases for $\sigma=0,0.01,0.05,$ and 0.10, but decreases from TPR = 0.4 to 0 at both $\sigma=0.25$ and 0.50. Thus, pruning always helped PDE-FIND learn the true equation when using the ANN method, and helps the other computational methods for small noise levels. 

For the nonlinear Fisher-KPP qquation, the median TPR value always improved the accuracy of the PDE-FIND method when using ANN approximations (Figure \ref{box_whisker_fisher_nonlin_no_prune}). When using finite difference approximations, the median TPR value increases for $\sigma=0,0.01.05$. The median value decreases from TPR = 0.3 to 0 at $\sigma=0.10$  for finite difference approximations. When using spline approximations, the median TPR value increases when $\sigma=0\text{ and }0.01$. The median TPR value decreased from TPR = 0.3 to 0 at $\sigma=0.50$. While the median value is never TPR = 1 for this equation, these results suggest that pruning in general helps reduce the number of incorrect terms in the library.

\section{Tables of learned PDEs}

This section contains tables of the final learned PDEs for data from each equation considered at a given noise level. The equation form is the one most commonly selected by the PDE-FIND method with pruning over the 1,000 different training-validation splits of $u_t$ and $\Theta$. The provided parameter values are the mean value for these parameters when the equation form was the final learned equation.

\begin{table}
    \centering
    \begin{tabular}{|c|c|c|}
    \hline 
     &  & \textbf{True Equation}\tabularnewline
    \hline 
     &  & $u_{t}=-0.8u_{x}+.01u_{xx}$\tabularnewline
    \hline 
    $\boldsymbol{\sigma}$ & \textbf{method} & \textbf{learned equation}\tabularnewline
    \hline 
    00 & FD & $u_{t}=-0.800u_{x}+0.010u_{xx}-0.000350u^{2}u_{x}$\tabularnewline
    \hline 
    01 & FD & $u_{t}=-0.800u_{x}+0.010u_{xx}$\tabularnewline
    \hline 
    05 & FD & $u_{t}=0$\tabularnewline
    \hline 
    10 & FD & $u_{t}=0$\tabularnewline
    \hline 
    25 & FD & $u_{t}=0$\tabularnewline
    \hline 
    50 & FD & $u_{t}=0$\tabularnewline
    \hline 
    00 & SP & $u_{t}=-0.808u_{x}+0.011u_{xx}+0.001u^{2}u_{x}+0.000u^{2}u_{xx}$\tabularnewline
    \hline 
    01 & SP & $u_{t}=-0.794u_{x}+0.012u_{xx}$\tabularnewline
    \hline 
    05 & SP & $u_{t}=-0.797u_{x}+0.012u_{xx}$\tabularnewline
    \hline 
    10 & SP & $u_{t}=-0.774u_{x}$\tabularnewline
    \hline 
    25 & SP & $u_{t}=-0.709u_{x}$\tabularnewline
    \hline 
    50 & SP & $u_{t}=0$\tabularnewline
    \hline 
    00 & ANN & $u_{t}=-0.809u_{x}+0.011u_{xx}$\tabularnewline
    \hline 
    01 & ANN & $u_{t}=-0.803u_{x}+0.011u_{xx}-0.000u^{2}u_{xx}+0.001u_{x}^{2}$\tabularnewline
    \hline 
    05 & ANN & $u_{t}=-0.810u_{x}+0.011u_{xx}$\tabularnewline
    \hline 
    10 & ANN & $u_{t}=-0.809u_{x}+0.010u_{xx}$\tabularnewline
    \hline 
    25 & ANN & $u_{t}=-0.796u_{x}+0.009u_{xx}$\tabularnewline
    \hline 
    50 & ANN & $u_{t}=-0.802u_{x}+0.008u_{xx}+0.001u_{x}^{2}$\tabularnewline
    \hline 
    \end{tabular}
    
    \caption{Learned Equations for the diffusion-advection equation.}
    \label{learned_eqns_diffadv}
\end{table}

\begin{table}
\centering
\begin{tabular}{|c|c|c|}
\hline 
 &  & \textbf{True Equation}\tabularnewline
\hline 
 &  & $u_{t}=0.02u_{xx}-10u^{2}+10u$\tabularnewline
\hline 
\textbf{$\boldsymbol{\sigma}$} & \textbf{method} & \textbf{learned equation}\tabularnewline
\hline 
00 & FD & $u_{t}=0.020u_{xx}-9.994u^{2}+9.996u$\tabularnewline
\hline 
01 & FD & $u_{t}=-10.155u^{2}+9.951u$\tabularnewline
\hline 
05 & FD & $u_{t}=0$\tabularnewline
\hline 
10 & FD & $u_{t}=0$\tabularnewline
\hline 
25 & FD & $u_{t}=0$\tabularnewline
\hline 
50 & FD & $u_{t}=0$\tabularnewline
\hline 
00 & SP & $u_{t}=0.020u_{xx}-9.993u^{2}+9.997u$\tabularnewline
\hline 
01 & SP & $u_{t}=0.020u_{xx}-9.972u^{2}+9.978u$\tabularnewline
\hline 
05 & SP & $u_{t}=-10.130u^{2}+9.926u$\tabularnewline
\hline 
10 & SP & $u_{t}=-10.088u^{2}+9.916u$\tabularnewline
\hline 
25 & SP & $u_{t}=0$\tabularnewline
\hline 
50 & SP & $u_{t}=0$\tabularnewline
\hline 
00 & ANN & $u_{t}=0.023u_{xx}-9.308u^{2}+9.533u$\tabularnewline
\hline 
01 & ANN & $u_{t}=0.020u_{xx}-9.972u^{2}+9.978u$\tabularnewline
\hline 
05 & ANN & $u_{t}=0.021u_{xx}-9.734u^{2}+9.837u$\tabularnewline
\hline 
10 & ANN & $u_{t}=0.022u_{xx}-9.287u^{2}+9.588u$\tabularnewline
\hline 
\multirow{2}{*}{25} & \multirow{2}{*}{ANN}& $u_{t}=0.012u_{xx}-11.161u^{2}+12.537u$\tabularnewline
 &  & $+0.071uu_{xx}-0.105u_{x}^{2}$\tabularnewline
\hline 
\multirow{2}{*}{50} & \multirow{2}{*}{ANN} & $u_{t}=-0.016u_{x}+0.014u_{xx}-8.689u^{2}+12.180u$\tabularnewline
 &  & $-0.034u^{2}u_{x}+0.077uu_{xx}-0.109u_{x}^{2}$\tabularnewline
\hline 
\end{tabular}

\caption{Discovered Equations for the Fisher-KPP Equation}
\label{learned_eqns_fisher_kpp}
\end{table}

\begin{table}
    \centering
    \begin{tabular}{|c|c|c|}
    \hline 
     &  & \textbf{True Equation}\tabularnewline
    \hline 
     &  & $u_{t}=-10u^{2}+10u+0.02uu_{xx}+0.02u_{x}^{2}$\tabularnewline
    \hline 
    \textbf{$\boldsymbol{\sigma}$} & \textbf{method} & \textbf{learned equation}\tabularnewline
    \hline 
    00 & FD & $u_{t}=0.000u_{xx}-9.996u^{2}+9.997u+0.020uu_{xx}+0.020u_{x}^{2}$\tabularnewline
    \hline 
    01 & FD & $u_{t}=-10.268u^{2}+10.211u$\tabularnewline
    \hline 
    05 & FD & $u_{t}=-9.532u^{2}+9.731u$\tabularnewline
    \hline 
    10 & FD & $u_{t}=0$\tabularnewline
    \hline 
    25 & FD & $u_{t}=0$\tabularnewline
    \hline 
    50 & FD & $u_{t}=0$\tabularnewline
    \hline 
    00 & SP & $u_{t}=0.001u_{xx}-10.013u^{2}+10.013u+0.019uu_{xx}+0.018u_{x}^{2}$\tabularnewline
    \hline 
    01 & SP & $u_{t}=0.006u_{xx}-9.821u^{2}+9.790u+0.018u_{x}^{2}$\tabularnewline
    \hline 
    05 & SP & $u_{t}=-10.393u^{2}+10.288u$\tabularnewline
    \hline 
    10 & SP & $u_{t}=-10.265u^{2}+10.195u$\tabularnewline
    \hline 
    25 & SP & $u_{t}=-10.146u^{2}+10.083u$\tabularnewline
    \hline 
    50 & SP & $u_{t}=0$\tabularnewline
    \hline 
    00 & ANN & $u_{t}=0.010u_{xx}-9.295u^{2}+9.238u-0.032uu_{xx}+0.017u_{x}^{2}$\tabularnewline
    \hline 
    01 & ANN & $u_{t}=-9.398u^{2}+9.390u+0.025u_{x}^{2}$\tabularnewline
    \hline 
    05 & ANN & $u_{t}=0.009u_{xx}-9.312u^{2}+9.246u-0.032uu_{xx}+0.017u_{x}^{2}$\tabularnewline
    \hline 
    10 & ANN & $u_{t}=-0.006u_{xx}-8.965u^{2}+9.140u+0.027u_{x}^{2}$\tabularnewline
    \hline 
    25 & ANN & $u_{t}=-0.010u_{xx}-7.928u^{2}+8.552u+0.034u_{x}^{2}$\tabularnewline
    \hline 
    50 & ANN & $u_{t}=0.286-0.026u_{xx}-5.419u^{2}+6.724u+0.045u_{x}^{2}$\tabularnewline
    \hline 
    \end{tabular}

\caption{Discovered Equations for the nonlinear Fisher-KPP Equation.} \label{learned_eqns_nonlinear_fisher_kpp}
\end{table}